\documentclass[12pt]{article} 
\usepackage{amsmath}
\newtheorem{theorem}{Theorem}[section]
\newtheorem{remark}[theorem]{Remark}
\newtheorem{definition}[theorem]{Definition}
\newtheorem{proposition}[theorem]{Proposition}
\newtheorem{example}[theorem]{Example}

\newtheorem{lemma}[theorem]{Lemma}

\begin{document}

\def\edit{\bf }
\def\endedit{\rm }
\def\square{ }

\def\Del{\Delta}
\def\L{\Lambda}
\def\Om{\Omega}
\def\S{\Sigma}
\def\s{\sigma}
\def\r{\rho}
\def\l{\lambda}
\def\b{\beta}
\def\a{\alpha}
\def\G{\Gamma}
\def\Ga{\Gamma}
\def\g{\gamma}
\def\k{\kappa}
\def\l{\lambda}
\def\ve{\varepsilon}
\def\v{\varepsilon}
\def\vp{\varphi}
\def\o{\omega}
\def\proof{\noindent {\bf Proof: \quad}}

\def\TCR{{\rm TCR}\,}
\def\CR{{\rm CR}\,}
\def\GaCR{{\rm \G CR}\,}
\def\dim{{\rm dim}\,}
\def\Ker{{\rm Ker}\,}
\def\soc{{\rm soc}\,}
\def\wr{{\rm wr}\,}

\def\Aut{{\rm Aut}\,}
\def\Inn{{\rm Inn}\,}
\def\Out{{\rm Out}\,}
\def\Sym{{\rm Sym}\,}
\def\ColPairs{{\rm ColPairs}}
\def\NonCol{{\rm NonCol}}
\def\DisjFlag{{\rm DisjointFlag}}

\def\AGL{{\rm AGL}\,}
\def\AG{{\rm AG}\,}
\def\PG{{\rm PG}\,}
\def\PSL{{\rm PSL}\,}
\def\PGL{{\rm PGL}\,}
\def\PSU{{\rm PSU}\,}
\def\PGL{{\rm PGL}\,}
\def\GF{{\rm GF}\,}
\def\GaL{{\rm \G L}\,}
\def\GaSp{{\rm \G Sp}\,}
\def\PGaL{{\rm P\G L}\,}
\def\GL{{\rm GL}\,}
\def\SL{{\rm SL}\,}
\def\Sp{{\rm Sp}\,}
\def\GSp{{\rm GSp}\,}
\def\GO{{\rm GO}}
\def\GU{{\rm GU}\,}
\def\U{{\rm U}\,}
\def\SU{{\rm SU}\,}
\def\O+{{\rm O}$^+$\,}
\def\O-{{\rm O}$^-$\,}
\def\Oe{{\rm O}$^\v$\,}

\def\qed{\hfill \rule{2mm}{2mm}}

\def\amalq{\stackrel{\mid}{\wedge}}

\def\aut#1{{\rm{ Aut}}(#1)} 
\def\diff{\mathbin{\mkern-1.5mu\setminus\mkern-1.5mu}}

\def\la{\langle}
\def\ra{\rangle}

\def\C{\mathcal{C}}
\def\calS{\mathcal{S}}
\def\B{\mathcal{B}}
\def\D{\mathcal{D}}
\def\K{\mathcal{K}}
\def\C{\mathcal{C}}
\def\H{\mathcal{H}}
\def\calP{\mathcal{P}}

\def\CC{{\Bbb C}}
\def\QQ{{\Bbb Q}}

\def\shat{{\hat{\s}}}

\title{Symmetric graphs with complete quotients} 

\date{Draft 5:\quad May 2017}

\date{Dedicated with affection and admiration to the memory\\ of our friend and colleague Anne Penfold Street}

\author{A. Gardiner\thanks{13, Bincleaves Road, Weymouth DT4 8RL, United Kingdom}
\and Cheryl E. Praeger\thanks{School of Mathematics and Statistics, 
University of Western Australia, Crawley, WA 6009, Australia. 
The research for this project was partially supported by an
Australian Research Council large grant. MSC codes: 20B25, 05C25.\qquad Draft 4.2:\quad April 2017}}
\maketitle

%
%
%


\bigskip\noindent 
{\bf Abstract}

\bigskip

Let $\G$ be a $G$-symmetric graph with vertex set
$V$. We suppose that $V$ admits a $G$-invariant partition $\B =\{ B=B_0,
B_1,\ldots ,B_b\}$, with parts
$B_i$ of size $v$, and that the quotient graph $\G_\B$ induced on 
$\B$ is a complete graph $K_{b+1}$. Then, for each pair of
suffices $i,j$ ($i\ne j$), the graph $\la B_i, B_j\ra$ induced on
$B_i\cup B_j$ is bipartite with each vertex of valency 0 or  $t$ (a constant). When
$t=1$, it was shown earlier how a flag-transitive 1-design $\D (B)$ induced
on the  part $B$ can sometimes be used to classify possible triples
$(\G , G, \B )$.  Here we extend these ideas to $t\ge 1$ and prove
that, if $G(B)^B$ is 2-transitive and the blocks of $\D (B)$ have
size less than $v$, then either (i) $v<b$, or 
(ii) the triple $(\G , G, \B )$ is known explicitly.

\section{Introduction}
\label{sec-intro}

A graph $\G$ is $G$-{\it symmetric} if $G\le\Aut\G$ acts
transitively on the vertices of $\G$, and the stabiliser $G(x)$ of
a vertex $x$ acts transitively on the edges incident with $x$. Let
$\G$ be a $G$-symmetric graph with vertex set $V$, and let  $\B =\{
B=B_0, B_1, B_2,\ldots \}$ be a $G$-invariant partition of $V$. Then we obtain a
natural quotient graph $\G_\B$ with vertex set $\B$ (where two parts
$B_i, B_j\in\B$ are adjacent in $\G_\B$ precisely when they are
joined by some edge of $\G$). In seeking to analyse $\G$ in terms of 
the $G$-invariant partition $\B$ and the quotient graph $\G_\B$, it is natural to
exclude the case where $\G_\B$ is a null graph. 

To investigate such triples $(\Ga, G, \B)$, a framework \cite{GP1} 
was introduced by the authors in 1995 which  made use also of two other associated 
combinatorial structures: namely a 1-design $\D(B)$ induced on a part $B\in\B$, 
and the (bipartite) restriction $\la B_i, B_j\ra$ of $\Ga$ to the 
union $B_i\cup B_j$ of adjacent parts $B_i, B_j$ of $\Ga_\B$. 
The subgraph $\la B_i, B_j\ra$ has edges involving $k$ vertices of
each of $B_i, B_j$, with $k$ independent of $i,j$. The paper 
\cite{GP1} studies varous properties of these graphs, and then focuses on the situations 
where $\Ga_\B$ is a complete graph or a cycle, and the cases where
$k$ is $1, 2$, $v-1$ or $v$, where $v=|B_i|$.  Many of these special cases 
(and also the case $k=v-2$ \cite{LPZ10, XZ})
are explored further in the literature, especially by Sanming Zhou and his co-authors.
In particular they classified all triples  $(\Ga, G, \B)$
such that $\Ga_\B$ is a complete graph and $k=v-1$, a culmination of work in
\cite{CZ, FFXZ, GMPZ, Z02}; and they studied more general quotients 
$\Ga_\B$, especially the case where $G$ acts 2-arc transitively on $\Ga_\B$, 
see  \cite{IPZ, LPZ, LZ, Z04, Z08, Z09}.

This paper was motivated by studies initiated in \cite{GP1}, and
developed further in \cite{GP2, cr}, of the case where $\Ga_\B$ 
is a complete graph $K_{b+1}$, and the induced action $G(B)^B$ of the stabiliser 
of a part $B\in \B$ on $B$ is 2-transitive. The results from these papers 
typically assert that either the triple  $(\Ga, G, \B)$ is known 
explicitly, or the part-size $|B|$
is bounded above by a certain function $f(b)$. For example, \cite[Theorem]{GP2}
proves this with $f(b)=b+1$ if the bipartite graph  $\la B_i, B_j\ra$  
is a perfect matching. An earlier special case (of \cite[Theorem]{GP2}),
namely \cite[Theorem 4.3]{GP1}, identified $G$-distance-transitive, 
antipodal covers of $K_{b+1}$ as examples (and all such covers are known,
see \cite{GLP, T}). 

Here we focus on the situation where
$\la B_i, B_j\ra$  has isolated vertices and may, or may not, be 
a partial matching. Our main result Theorem~\ref{main} 
identifes a rich class of examples in this case, and finds all examples
for which $|B|\geq b$. 
Two of the families of examples (first identified in an early version of this 
paper \cite{GP3}), namely the cross ratio graphs and twisted cross ratio graphs,
have been studied in detail in \cite{cr}; their full automorphism 
groups are determined there, as are the small number of exceptional isomorphisms 
between them. First we summarise briefly some fundamental theory of $G$-imprimitive symmetric graphs, 
and then we state Theorem~\ref{main}.

\subsection{Basic facts and parameters}\label{sub:basic}

Let $(\Ga, G, \B)$ be as above, that is, $\Ga = (V, E)$ is a $G$-symmetric 
(simple undirected) graph with vertex set $V$ and edge set $E$, 
$\B$ is a $G$-invariant vertex-partition, and the quotient $\Ga_\B$ is 
not a null graph. Since $G$ is transitive on $V$, $\G$ is regular of valency
$s=|\Ga(x)|$, where $\Ga(x)=\{y | \{x,y\}\in E\}$; 
all parts $B_i\in\B$ have the same size $v:=|B|$; and the
setwise stabiliser $G(B)$ of $B\in\B$ acts transitively on $B$. Since
$G$ acts transitively on the edges of $\G$ and since $\G_\B$ is not 
a null graph, the number of edges joining
each pair of adjacent parts $B_i, B_j$ is a non-zero constant $m$ (say),
and there are no edges joining two vertices in the same part. If $K$
is the kernel of the action of $G$ on $\B$, then the quotient
graph $\G_\B$ is $G/K$-symmetric, so
$\G_\B$ is regular of valency $b$, where $v\cdot s=b\cdot m$. 

If $B_i, B_j$ are adjacent parts, then since $G$ acts transitively
on the \emph{arcs} of $\G$ (ordered pairs of adjacent vertices), 
each vertex of $B_i$ is adjacent to either $0$ or
$t$ vertices of $B_j$, for some constant $t\ge 1$. Thus $t$ divides
both $s$ and $m$. Set 
\begin{equation}\label{eq-1}
r:=s/t\quad \hbox{and}\quad k:=m/t. 
\end{equation} 
The study in~\cite{GP1} looked solely at the case where $t=1$, that is, where  
$\la B_i, B_j\ra$ ($i\ne j$) is a partial
matching. However, in general, the parameter $t$ can be arbitrarily large (even
if we demand not only that $\G$ is $G$-symmetric, but also that
$G(B)^B$ is 2-transitive and that $\G_\B$ is complete: see 
Example~\ref{ex-2.8}, and line 3 of Table~\ref{tab-main1}). 

For $B\in\B$, set $\Ga(B):=\{ y\, |\, \Ga(y)\cap B\ne\emptyset\}$
and $\Ga_\B(B):=\{ B'\in\B \,|\, \Ga(B')\cap B\ne\emptyset\}$. 
The $1$-design $\D(B)$ defined in \cite{GP1} has point set $B$,
block set $\Ga_\B(B)$, and a `block' $B'$ is adjacent to a `point' $x$ if
$\Ga(x)\cap B'\ne\emptyset$. We define the \emph{parameters} of $(\Ga,G,\B)$ to be
$(v,b,r,k,s,t)$ which we summarise in Table~\ref{tab-param}.

\begin{table}
\begin{center}
 \begin{tabular}{|llll|}\hline
  $v$	&$=$	&	$|B|$	& the number of points of $\D(B)$\\ 
  $b$	&$=$	&	$|\Ga_\B(B)|$	& the valency of $\Ga_B$ and the number of blocks of $\D(B)$\\ 
  $r$	&$=$	&	$s/t$	& the number of blocks on a point of $\D(B)$\\ 
  $k$	&$=$	&	$|\Ga(B')\cap B|$& if $y\in B'$ and  $\Ga(y)\cap B\ne \emptyset$; the block-size of $\D(B)$\\ 
  $s$	&$=$	&	$|\Ga(x)|$	& the valency of $\Ga$\\ 
  $t$	&$=$	&	$|\Ga(x)\cap B'|$& if $x\in B, \Ga(x)\cap B'\ne \emptyset$\\ \hline
 \end{tabular}
\caption{Parameters of $(\Ga, G,\B)$}\label{tab-param}
\end{center}
\end{table}


\subsection{Main result}

Our goal, which we partially achieve in this paper, is as follows: 

\begin{quote}
\emph{For arbitrary $t\ge 1$ and arbitrary $k\le v$, we would like to classify all
`exceptional' triples $(\G , G, \B )$ with $G(B)^B$ $2$-transitive and  $\G_\B=K_{b+1}$,
and with $v\ge b$.}  
\end{quote}
In~\cite{GP2} we
explained why such a result is of interest (for example, as an
extension of Fisher's inequality $v\le b$ for 2-designs, and as a generalisation of the
inequality $r\le k$ for an $r$-fold antipodal distance transitive
cover of a $k$-valent graph). In this paper we accomplish this when $k < v$.
The case $k=v$ remains open.

First we comment on the notation used in the statement of Theorem~\ref{main}.
\begin{enumerate}
 \item[$\bullet$]  We denote by $c\cdot \Del$
the graph consisting of $c$ disjoint copies of a given graph $\Del$.

\item[$\bullet$] The complete graph on $n$ vertices is denoted $K_n$ while $K_{n[a]}$
denotes the complete multipartite graph with $n$ parts of size $a$;
if $n=2$ we usually write $K_{2[a]}$ as $K_{a,a}$ (a complete bipartite graph). 

\item[$\bullet$] Several of the examples 
listed are ``$*$-transforms'' of more familiar graphs (see
Definition~\ref{def-2.6}).

\item[$\bullet$] The cross-ratio graphs $\CR(q;d,s)$ and $\TCR(q;d,s)$ 
arising in Theorem~\ref{main}~(c) are defined in Definition~\ref{def-2.1}. 
For these graphs $d\in\GF(q)\setminus\{0,1\}$, and $s$ divides $s(d)$,
where the subfield of $\GF(q)$ generated by $d$ has order $p^{s(d)}$
(with $p$ the prime dividing $q$).  The relevant groups $G$ 
for these graphs are subgroups of $\PGaL(2,q)$ which are $3$-transitive 
on $\GF(q)\cup\{\infty\}$; such subgroups are classified in \cite[Theorem 2.1]{cr}.

\item[$\bullet$] The graphs associated with designs, which occur 
in the tables for Theorem~\ref{main} are introduced in Subsection~\ref{sub-designs}.
In particular, the collinear pairs graph $\ColPairs(\D')$, where 
$\D'$ is the design of points and planes in the affine space $\AG(d,2)$,
is isomorphic to $(2^d-1)\cdot K_{2^{d-1}[2]}$. The design $S(22,6,1)$ on line 2 
of Table~\ref{tab-main2} is the Steiner system with automorphism group $\Aut(M_{22})$,
and the design $\D(M_{11})$ on line 3 of Table~\ref{tab-main2}
is the unique $3-(12,6,2)$ design with automorphism group $M_{11}$.

\end{enumerate}

\begin{theorem}\label{main}%
Let $\G=(V,E)$ be a $G$-symmetric graph with a $G$-invariant 
vertex-partition $\B =\{B_0, B_1,\ldots ,B_b\}$ and parameters 
$(v,b,r,k,s,t)$. Suppose that the quotient $\G_\B$ is a complete
graph $K_{b+1}$, and that, for $B\in\B$, $G(B)^B$ is $2$-transitive, 
and $k<v$. Then one of the following holds.
\begin{description}
 \item[(a)] $v<b$; 

 \item[(b)]$v=b$, $r=k$, $V=\{ ij\ :\ i, j\in X,\ i\ne j\}$ where 
 $X=\{0,1,\dots,v\}$, $B_i=\{ ij\
:\ j\in X, j\ne i\}$ (for $i\in X$), $G\leq \Sym(X)$ acts 
coordinate-wise,  and $\Ga, G, t$ are as in one of the lines of 
Table~$\ref{tab-main1}$, where $k=1$ in line $1$, and $k=v-1$ in the other lines. 

\item[(c)] $v=b$, $V$ is the set of flags $i\b$ of a design $\D'$ with point set $X=\{0,1,\dots,v\}$,
$B_i$ is the set of flags $i\b$ on $i$ (for $i\in X$), $G\leq \Aut(\D')\leq \Sym(X)$
acts coordinate-wise, and $\Ga, G, \D', t$ are as in one of the lines of 
Table~$\ref{tab-main3}$.
  \end{description}
\end{theorem}

\begin{table}
 \begin{center}
\begin{tabular}{lllll}
Line & $\Ga$      &Edges $\{ij,i'j'\}$  	&$t$& Conditions on $G$ and parameters  \\ \hline
1    & $\binom{v+1}{2}\cdot K_2$&$(i',j')=(j,i)$&$1$& $G$ $3$-transitive on $X$ \\
2    & $(v+1)\cdot K_v$         &$j=j'$         &$1$& $G$ $3$-transitive on $X$\\
3    & $*$-transf. of line 2  &Def.~\ref{def-2.6}&$v-2$& $G=A_{v+1}$ ($v\geq 5$), $S_{v+1}$ ($v\geq4$) \\
  &				&		&	& or $M_{v+1}$ ($v+1=11, 12, 23, 24$) \\
4    & $\G = \CR (v;d,s(d))$    &Def.~\ref{def-2.1}&$s/s(d)$& $r=k=v-1$, $v$ a prime power,\\ 
  &	or $\TCR(v;d,s(d))$	&		&	& $G\leq \PGaL(2,v)$, $3$-trans. on $X$, \\
  &				&		&	& $G(\infty 01)$-orbit cont'g $d$ has size $t$ \\
5    & $\ColPairs(\D')$     &Def~\ref{def-desgrph}&$t_1$&  $\D', G, t_1$ as in Table~\ref{tab-main2}, lines $1, 2, 3$\\
6    & $\NonCol(\D')$     &Def~\ref{def-desgrph}&$t_2$&  $\D', G, t_2$ as in Table~\ref{tab-main2}, lines $1, 2, 3$\\ \hline
\end{tabular}
  \end{center}
  \caption{Examples for Theorem~\ref{main}~(b)}\label{tab-main1}
\end{table}

\begin{table}
 \begin{center}
\begin{tabular}{lllll}
Line & $\Ga$      &Edges $\{i\beta,i'\beta'\}$  	&$t$& Conditions on $G$ and parameters  \\ \hline
1    & $b'\cdot K_{k'}$&$\beta=\beta'$ &$1$&  $\D', G, b', k'$ as in Table~\ref{tab-main2}, lines $2$--$4$ \\
2    & $(b'/2)\cdot K_{k',k'}$&$\b\cap\b'=\emptyset$ &$1$& $\D', G, b', k'$ as in Table~\ref{tab-main2}, lines $3, 4$\\
3    & $*$-transf. of line 1  &Def.~\ref{def-2.6}&$k'-2$& $\D', G, b', k'$ as in Table~\ref{tab-main2}, lines $2$--$4$ \\
4    & $*$-transf. of line 2  &Def.~\ref{def-2.6}&$k'-1$& $\D', G, b', k'$ as in Table~\ref{tab-main2}, lines $3, 4$ \\
5    & $\Ga_1(M_{22})$  &$i\not\in\b', i'\not\in \b$, &$6$& $\D', G$ as in Table~\ref{tab-main2}, line $2$ \\
      &			&$\b\cap\b'=\emptyset$	      &	  &  \\
6    & $\Ga_2(M_{22})$  &$i\not\in\b', i'\not\in \b$, &$10$& $\D', G$ as in Table~\ref{tab-main2}, line $2$ \\
 &			&$|\b\cap\b'|=2$	&	&  \\ \hline
\end{tabular}
  \end{center}
  \caption{Examples for Theorem~\ref{main}~(c); in lines $3$--$6$, $k'=k+1$}\label{tab-main3}
\end{table}

\begin{table}
 \begin{center}
\begin{tabular}{lllllllll}
Line&$\D'$		& $v+1$		& $k'$	& $\lambda$& $G$	& $t_1$& $t_2$ & $b'$\\ \hline
$1$&${\rm AG}_2(d,2)$& $2^d$	& $4$	& $1$	& $\AGL(d,2)$ , or 	& $1$& $2^d-4$&\\
&		&		&	&	&$d=4, G=Z_2^4\cdot A_7$&  & &\\
$2$&$S(22,6,1)$	&$22$		& $6$	& $1$	& $M_{22}$ or $\Aut(M_{22})$& $3$&$16$& $77$\\
$3$&$\D'(M_{11})$	&$12$		& $6$	& $2$	& $M_{11}$		& $6$&$3$&$22$\\ 
$4$& ${\rm AG}_{d-1}(d,2)$& $2^d$	& $2^{d-1}$&$2^{d-2}-1$	& Same as line $1$	& --&-- &$2^{d+1}-2$\\ \hline 
\end{tabular}
  \end{center}
  \caption{Designs $\D'$ for Theorem~\ref{main}~(b) and~(c); $\D'$ has point set $X$ and $b'$ blocks
  and is a $3-(v+1,k',\lambda)$ design}\label{tab-main2}
\end{table}

In Section~\ref{sec-2} we introduce the main families of triples $(\Ga, G, \B)$
that arise in our classification. In Section~\ref{sec-3} we 
indicate how the framework which was introduced in~\cite{GP1}, 
for the case where $t=1$ and $G(x)^{\G (x)}$ is primitive, can
be extended to arbitrary symmetric graphs with $t\ge 1$. In particular
we show that in any such graph we get a 1-design $\D (B)$ induced on
the part $B$, with $G(B)$ acting flag-transitively on $\D (B)$.
This framework is then used in Sections~\ref{sec-4} and~\ref{sec-5} 
to analyse triples 
$(\G , G, \B )$ in which $G(B)^B$ is 2-transitive and $\G_\B
=K_{b+1}$ is a complete graph. We treat separately the cases 
$t=1$ (Section~\ref{sec-4}) and $t\ge 2$
(Section~\ref{sec-5}).

\section{Examples and general constructions.}\label{sec-2}%

In this section we introduce the main families of graphs that arise in the 
situation analysed in Theorem~\ref{main}. We
begin by defining the cross-ratio 
graphs. Next we look at some generic examples arising from the action
of a 3-transitive permutation group $(G,\Om)$ on the set of ordered pairs 
of distinct elements of $\Om$, and from the flags in a 3-design.
We then introduce and illustrate the ``$*$-transform'' construction
which links certain pairs of triples $(\G, G, \B)$ and $(\G^*,G,\B)$. 

\subsection{Cross ratio graphs}\label{sub-cr}

First we introduce the {\it untwisted cross ratio graphs} $\CR(q;d,s)$ and the
{\it twisted cross ratio graphs} $\TCR(q;d,s)$ which feature in
Sections~\ref{sec-4} and~\ref{sec-5}.  They may be defined as orbital graphs of
a transitive permutation group with respect to nontrivial self-paired orbits of
a point stabiliser.  If $G$ is a transitive permutation group on a set $V$, and
$x\in V$, then the $G(x)$-orbit $Y$ containing a point $y$ is {\it nontrivial}
provided $y\ne x$, and is {\it self-paired} if there exists an element in $G$
which interchanges $x$ and $y$.  For a nontrivial self-paired $G(x)$-orbit $Y$,
the corresponding {\it orbital graph} is defined as the graph with vertex set
$V$ and edges the pairs $\{ x^g, y^g\}$, for $y\in Y$ and $g\in G$.  It is
straightforward to show that each orbital graph for $G$ is $G$-symmetric.

Let $q=p^n$ where $p$ is prime and $n$ is a positive integer.  The projective
line $\PG(1,q)$ over the field $\GF(q)$ of order $q$ can be identified with the
set $\GF(q)\cup\{\infty\}$, where $\infty$ satisfies the usual arithmetic rules
such as $1/\infty = 0, \infty + y = \infty$, etc.  The two-dimensional
projective group $\PGL(2,q)$ then consists of all fractional linear
transformations
\[
t_{a,b,c,d}: z \mapsto \frac{az+b}{cz+d}\,  \quad (\mbox{with} \ 
a,b,c,d \in \GF(q), \ \mbox{and}\ ad - bc \neq 0)
\]
of $\PG(1,q)$ (see, for example, \cite[p. 242]{dm}). Note that
$t_{a,b,c,d}=t_{a',b',c',d'}$ if and only if the 4-tuple $(a,b,c,d)$ is a non-zero
multiple of $(a',b',c',d')$. The group $\PGL(2,q)$ is sharply 3-transitive
in this action on $\PG(1,q)$, that is, it is 3-transitive and only the identity
element $t_{1,0,0,1}$ fixes three elements of $\PG(1,q)$.
The Frobenius automorphism $\s:x\mapsto x^p$ of the field $\GF(q)$ induces an
automorphism of $\PGL(2,q)$ by $\s:t_{a,b,c,d}\mapsto t_{a^p,b^p,c^p,d^p}$, and
the group generated by $\PGL(2,q)$ and $\s$ is the semidirect product
$\PGL(2,q)\cdot\la\s\ra$ and is denoted by $\PGaL(2,q)$.  
The group $\PGaL(2,q)$ is the automorphism
group of $\PGL(2,q)$, and it too acts on $\PG(1,q)$ (with $\s:z \mapsto z^p$,
where $\infty^p=\infty$). The 3-transitive subgroups of $\PGaL(2,q)$ are
(see~\cite[Theorem 2.1]{cr}) precisely the subgroups $\PGL(2,q)\cdot\la\s^s\ra$
for divisors $s$ of $n$, and, in the case where $p$ is odd and $n$ is even,
also the subgroups ${\rm M}(s,q)=\la\PSL(2,q),\s^st_{a,0,0,1}\ra$, for
divisors $s$ of $n/2$, where $a$ is a primitive element of $\GF(q)$.

Each 3-transitive subgroup of $\PGaL(2,q)$ acts transitively on the set $V$ of
ordered pairs of distinct points of the projective line $\{\infty\}\cup\GF
(q)$.  The cross-ratio graphs are certain orbital graphs for the actions on $V$
of  3-transitive subgroups of $\PGaL(2,q)$. 
The stabiliser in $\PGL(2,q)$ of $\infty 0$ is cyclic of order $q-1$ and
consists of the transformations $z\mapsto az$ ($a\in\GF (q)^\#$); while if $q$
is odd then the stabiliser in $\PSL(2,q)$ of $\infty 0$ consists of the
transformations $z\mapsto az$ with $a$ a square.  Let $G$ be a 3-transitive
subgroup of $\PGaL(2,q)$: namely either $G= \PGL(2,q)\cdot\la\s^s\ra$ for
some divisor $s$ of $n$, or $p$ is odd, $n$ is even, and $G={\rm M}(s/2,q)$
for some even divisor $s$ of $n$. Then the stabiliser $G(\infty 0)$ 
is transitive on $\PG(1,q)\setminus \{\infty,0\}$, and $G(\infty 0 1)$ is 
$\la\s^s\ra$ (see~\cite[Corollary 2.2]{cr}).

Let $q\geq 3$.  Then each element $d\in\GF(q)\setminus\{0,1\}$ generates a
subfield of $\GF(q)$ of order $p^{s(d)}$ for some divisor $s(d)$ of $n$.
Suppose that $d$ is such that $s$ divides $s(d)$, so that in particular $s(d)$
is even if $G={\rm M}(s/2,q)$.  Then the $G(\infty 0 1)$-orbit in $V$
containing $1d$ is the set $\{1d^{\s^{si}}\ \mid\ 0\leq i<s(d)/s\}$ of size
$s(d)/s$.  If $G= \PGL(2,q)\cdot\la\s^s\ra$, then the $G(\infty 0)$-orbit
$\Del(\infty 0)$ containing $1d$ consists of the $(q-1)s(d)/s$ pairs $ab$,
where $a\in\GF(q) \setminus\{0\}$ and $b=ad^{\s^{si}}$ with $0\leq i<s(d)/s$.
The orbit $\Del(\infty 0)$ is self-paired because the element $t_{1,-d,1,-1}$
interchanges $\infty 0$ and $1d$.  If $G= {\rm M}(s/2,q)$, then the $G(\infty
0)$-orbit $\Del'(\infty 0)$ containing $1d$ consists of the $(q-1)s(d)/s$ pairs
$ab$, where $b=ad^{\s^{si}}$ with $0\leq i<s(d)/s$ if $a$ is a square in
$\GF(q)$, and $b=ad^{\s^{si+s/2}}$ with $0\leq i<s(d)/s$ if $a$ is a non-square
in $\GF(q)$.  If $d-1$ is a square then the orbit $\Del'(\infty 0)$ is
self-paired because the element $t_{1,-d,1,-1}\in \PSL(2,q)\subseteq G$.
However if $d-1$ is a non-square then $\Del'(\infty 0)$ is not self-paired.
(This may be proved by an argument similar to that used in the last paragraph
of the proof of \cite[Theorem 4.1]{cr}.)  We are now able to define the
cross-ratio graphs.

\begin{definition}\label{def-2.1}%
{\rm 
Let $q=p^n$ for a prime $p$, where $n\geq 1$ and $q\geq 3$, and let $V$ denote
the set of ordered pairs of distinct points from the projective line $\PG(1,q)
=\GF(q)\cup\{\infty\}$. 
Let $d\in\GF(q)$, $d\ne 0, 1$, and let $s$ be a divisor of $s(d)$. 
\begin{enumerate}
\item[(a)] The \emph{untwisted cross ratio graph} $\CR(q; d,s)$ is 
defined as the orbital graph for the action on $V$ of $G:=\PGL(2,q)\cdot\la\s^s\ra$
corresponding to the self-paired $G(\infty 0)$-orbit $\Del(\infty 0)$.
\item[(b)] Suppose now that $p$ is odd, $n$ is even (so $q\geq 9$), 
and also that $d-1$ is a 
square and both $s$ and $s(d)$ are even. Then the \emph{twisted cross ratio 
graph} $\TCR(q; d,s)$ is defined as the orbital graph for the action on $V$
of the group $G:={\rm M}(s/2,q)$ corresponding to the self-paired
$G(\infty 0)$-orbit $\Del'(\infty 0)$.
\end{enumerate}
}
\end{definition} 

\begin{remark}\label{rem-cr}%
{\rm 
(a) The graphs $\CR(q;d,s)$ are defined in \cite[Definition 3.2]{cr} in terms
of the cross-ratio. However their definition given here as orbital graphs
is equivalent, see \cite[Remark 3.3 (b)]{cr}.

(b) The graphs $\CR(3;2,1)$ and $\CR(5;4,1)$ are disconnected: $\CR(3;2,1)
\cong 3\cdot C_4$ and $\CR(5;4,1)\cong 5\cdot (C_3[\overline{K_2}])$ 
(a lexicographic product). All other twisted and untwisted cross-ratio 
graphs are connected, by \cite[Proposition 5.2]{cr}.

(c)  All isomorphisms 
between cross-ratio graphs are specified in \cite[Theorem 6.2]{cr}. In particular
there are no isomorphisms between a twisted cross-ratio graph and an 
untwisted cross-ratio graph.  

(d) The full automorphism groups of all twisted and untwisted cross-ratio graphs
are determined in \cite[Theorem 6.1]{cr}. For $\G=\CR(q;d,s)$ we have $\Aut(\G)\cap\PGaL(2,q)=\PGL(2,q)\cdot
\la\s^s\ra$ by \cite[Theorem 3.4]{cr}, and for $\G=\TCR(q;d,s)$ we have 
$\Aut(\G)\cap\PGaL(2,q)={\rm M}(s/2,q)$ by \cite[Theorem 3.7]{cr}. In either case,
$G:=\Aut(\G)\cap\PGaL(2,q)$ acts symmetrically on $\G$ and preserves the 
$G$-invariant partition 
\[
\B:=\{  B(x)\ \mid \ x\in\PG(1,q)\}, \ \mbox{where}\  
B(x)=\{ xy\ \mid\ y\in\PG(1,q), y\ne x\}.
\]
The part size is $v=|B(x)|=q\geq 3$ and the quotient $\G_\B\cong K_{q+1}$.
Moreover, by \cite[Theorems 3.4 and 3.7]{cr}, the parameters $k=q-1$ and $t
= s(d)/s$.
}
\end{remark}

\subsection{Generic examples from a $3$-transitive group}\label{sub-designs}

We now give some examples which are in many ways typical of 
the triples $(\G,G,\B)$ which we shall meet in later sections. 
The examples given reflect the fact that we are primarily
interested in cases where $G(B)^B$ is 2-transitive.
Our first two examples are of triples $(\G, G, \B)$ in which $G^\B$
is 3-transitive, and in which the vertices of $\G$ are labelled by
ordered pairs of distinct parts of $\B$. Such a situation is
typical of the case where $k=v-1$ (Proposition~\ref{prop-4.1} 
and~\ref{prop-5.1}). The first example, which extends naturally 
to the case where the set $\calP$ of labels is the set of points of 
the affine geometry $\AG(d,2), d\geq 4$, illustrates the case $t=1$
(as in Section~\ref{sec-4}); the second and third examples illustrate 
the case $t\geq 2$ (as in Section~\ref{sec-5}).

\begin{example}\label{ex-2.3}%
{\rm
Let $\calP =\{ 0,1,\ldots,7\}$ be the set of points of the affine 
geometry $\AG(3,2)$, and let $V=\{ ij :\ 0\le i,j\le 7, i\ne j\}$. 
Define a graph $\G$ on the vertex set $V$, by joining vertex $ij$ to 
vertex $kl$ precisely when $i,j,k,l$ are the four points of 
some $2$-dimensional affine subspace. Then $G=\AGL(3,2)$ acts
$3$-transitively on $\calP$, $\G$ is $G$-symmetric, and we have a
natural $G$-invariant partition $\B =\{ B=B_0, B_1,\ldots,B_7\}$ on $V$, 
where $B_i = \{ ij\ :\ 0\le j\le 7, j\ne i\}$. Moreover $G(B)^B$
is $2$-transitive. The pair $(\G,\B)$ has parameters $v=b=7, 
r=k=6, t=1$. The graph $\G$ is the disjoint union of $7$ copies 
of the complete multipartite graph $K_{4;2}$ with $4$ parts of size $2$, 
and the quotient $\G_\B = K_8$.
}
\end{example}

\begin{example}\label{ex-2.4}%
{\rm
The group $G=S_5$ acts naturally on the set $\{ 0,1,2,3,4\}$.  Let the graph
$\G$ have vertex set $V=\{ ij :\ 0\le i,j\le 4, i\ne j\}$, with vertex $ij$
adjacent to vertex $kl$ if and only if $|\{ i,j,k,l\}|=4$.  The set $V$ admits
the $G$-partition $\B =\{ B=B_0, B_1,B_2, B_3, B_4\}$, where $B_i = \{ ij\ :\
0\le j\le 4, j\ne i\}$ with $G(B)^B = S_4$.  This defines a connected
$G$-symmetric graph of valency $6$, and a triple $(\G, G, \B)$ with parameters
$v=b=4, r=k=3, t=2$, and with quotient graph $\G_\B = K_5$.

In fact $\G = \bar{K_2} \wr O_3$ is obtained from the Petersen graph $O_3$ by
replacing each vertex $v$ of $O_3$ by two vertices $v_1, v_2$, and each edge
$vw$ of $O_3$ by the four edges $v_iw_j$ ($1\le i,j\le 2$).  The full
automorphism group $\Aut \G = S_2\wr S_5 = S_2^5\cdot S_5$ admits only the
invariant partition $\{ \{ v_1,v_2\}\ :\ v\in VO_3\}$. However the subgroup 
$G:= \langle A_5, d(12)\rangle \cong S_5$ (where $d$ is the generator of the 
diagonal subgroup of $S_2^5$ and $(12)$ lies in the top group $S_5$)  acts 
transitively on $1$-arcs, and the subgroup chain $G >
G(B) > G(v_1)$ (with both containments proper) gives rise to the $G$-invariant partition
$\B$ (compare~\cite[Example~2.3]{GP1}).  
}
\end{example}

\begin{example}\label{ex-2.5}%
{\rm
Let $V$ be the set of flags, that is, the incident point-hyperplane
pairs $i\b$, in the affine geometry $\AG(3,2)$.
Define a graph $\G$ on $V$ by joining $i\b$ to $i'\b'$ 
if and only if $i\ne i', \b\ne\b', i\in\b', i'\in\b$.
Then $\G$ is a connected graph of valency $6$. The group $G=\AGL(3,2)$ 
acts symmetrically on $\G$, and preserves the partition $\B$ 
with $8$ parts $B_i$ where, for any given point $i$, 
$B_i =\{ i\b\ :\ \b$ a hyperplane 
containing $i\}$ has size $v=7$, $k=3$, and $t=2$. The quotient graph
$\G_\B = K_8$.
}
\end{example}

The examples arising in Theorem~\ref{main} are often associated 
with designs. For a positive integer $s$, an \emph{$s$-$(v', k', \lambda)$ design} $\D'=(X,\mathcal{L})$,
where $s\leq k'\leq v'$,
consists of a set $X$ of cardinality $v'$, whose 
elements are called points, and a set 
$\mathcal L$ of $k'$-element subsets of $X$, called blocks, 
such that each $s$-element subset of $X$ is contained in exactly 
$\lambda$ blocks. The relevant designs $\D'$ for the statement of Theorem~\ref{main}
are $3$-designs and are listed in 
Table~\ref{tab-main2}; in these cases the cardinality of $\mathcal{L}$ is
$b'= \lambda \frac{v'(v'-1)(v'-2)}{k'(k'-1)(k'-2)}$. 
We will meet the following graphs associated with designs.

\begin{definition}\label{def-desgrph}
{\rm
Let $\D'=(X,\mathcal{L})$ be a $3$-$(v', k', \lambda)$ design with
$b'$ blocks. Let $V=\{ ij\ :\ i, j\in X,\ i\ne j\}$, and let 
$F=\{ i\b\ :\ i\in X,\ \b\in\mathcal{L},\ i\in \b\}$ (the set of 
flags of $\D'$).
\begin{description}
 \item[(a)] The \emph{collinear pairs graph} $\ColPairs(\D')$ has vertex set $V$, 
 and vertices $ij$ and  $i'j'$ are adjacent if and only if $i, j, i', j'$ 
 are pairwise distinct and are all contained in some block in $\mathcal{L}$. 

 \item[(b)] The \emph{non-collinear pairs graph} $\NonCol(\D')$ has vertex set $V$, 
 and vertices $ij$ and  $i'j'$ are adjacent if and only if $i, j, i', j'$ 
 are pairwise distinct but no block of $\mathcal{L}$ contains all of them. 
 
 \item[(c)] Let $\D'$ be the $3$-$(22, 6, 1)$ Steiner system in line 2 
 of Table~\ref{tab-main2}, so $\Aut(\D')=\Aut(M_{22})$. The graphs 
 $\Ga_1(M_{22})$ and $\Ga_2(M_{22})$ in lines 5 and 6 of Table~\ref{tab-main3}
 both have vertex set $F$, and two flags $i\b$ and $i'\b'$ are adjacent  
 if and only if $i\not\in\b', i'\not\in\b$, and either $\b\cap\b'=\emptyset$
 (in $\Ga_1(M_{22})$), or $|\b\cap\b'|=2$ (in $\Ga_2(M_{22})$). 
\end{description}

} 
\end{definition}

\subsection{The star-transform}

Next we introduce the ``$*$-transform'' of a triple $(\G, G, \B)$
which relates certain triples with $t\geq 2$ in Theorem~\ref{main}
to triples with $t=1$.
Given any triple $(\G, G, \B)$, if $B_i,
B_j\in\B$ are adjacent in $\Ga_\B$, then $|\G (B_i)\cap B_j|=k$, 
with $k$ as in equation~(\ref{eq-1}), see Table~\ref{tab-param}. 

\begin{definition}\label{def-2.6}%
{\rm
Let $\Ga$ be a finite graph and let 
$\B =\{ B_0, B_1, B_2, \ldots\}$ be a partition of the vertex set $V$. 
For $ i\ne j$, let $X_{ij} =\G(B_i)\cap B_j$, and suppose that,
whenever $X_{ij}\ne\emptyset$ the cardinality $|X_{ij}|\geq2$.
Then the \emph{$*$-transform} $\Ga^{*}$ of $\Ga$ relative to $\B$ 
has vertex set $V$, and edges all pairs $\{x,y\}$ such that,
for some $i\ne j$, $x\in X_{ij}, y\in X_{ji}$, and
$\{x,y\}$ is not an edge of $\Ga$.
Note that $(\G^*)^*=\G$.
}
\end{definition}

Thus, if $A(\Ga)$ denotes the set of arcs of $\Ga$, then 
$A(\Ga^*)\cap (B_i\times B_j) = (X_{ji}\times X_{ij})
\setminus A(\Ga)$.

\begin{lemma}
Let $\G$ be a finite $G$-symmetric graph admitting a 
$G$-invariant partition $\B$ 
of the vertex set $V$, where $\B =\{ B_0, B_1, B_2, \ldots\}$,
such that $(\Ga,G, \B)$ has parameters $(v,b,k,r,s,t)$, with $k\geq2$.
Let $\Ga^{*}, X_{ij}$ be as in Definition~$\ref{def-2.6}$, 
relative to $\B$. 
\begin{enumerate}
 \item[(a)] Then $G\leq \Aut(\Ga^*)$ and, whenever $B_i, B_j$ 
 are adjacent in $\Ga_\B$, $G(B_i,B_j)$ is transitive on 
 $A(\Ga)\cap (B_i\times B_j) = A(\Ga)\cap (X_{ji}\times X_{ij})$.

 \item[(b)] $\Ga^*$ is $G$-symmetric if
and only if, 
for some $i, j$ such that $X_{ij}\ne\emptyset$,
the stabiliser $G(B_i, B_j)$ is transitive on 
$(X_{ji}\times X_{ij})\setminus A(\Ga)$.

\item[(c)] If $t=k$ then $\Ga^*$ is a null graph, while if $t<k$
and $\Ga^*$ is $G$-symmetric, then 
$\Ga^*$ has parameters $(v^*, b^*, r^*, k^*, s^*, t^*)
=(v,b,r,k,rk-s,k-t)$.
 \end{enumerate}
\end{lemma}

\proof
It is straightforward to show that $G\leq \Aut(\Ga^*)$.
From the definition of the $X_{ij}$ it follows that 
$A(\Ga)\cap (B_i\times B_j) = A(\Ga)\cap (X_{ji}\times X_{ij})$,
the set of arcs $(x,y)$ with $x\in B_i$ and $y\in B_j$.
If $(x,y), (x',y')$ are two arcs in this set, then since $G$ 
is transitive on $A(\Ga)$, there is an element $g\in G$ such 
that $(x',y')=(x,y)^g=(x^g,y^g)$. Since $\B$ is $G$-invariant, 
and since $x, x'\in B_i$ we have $B_i^g=B_i$, and 
similarly $B_j^g=B_j$. Thus $g\in G(B_i,B_j)$, proving
part (a). 

If $\Ga^*$ is a null graph, then it is $G$-symmetric since 
$G$ is transitive on vertices, and also the local transitivity 
property holds vacuously, so the equivalence holds in this case.
Suppose now that $\Ga^*$ is not null, and 
suppose first that $\Ga^*$ is $G$-symmetric. Then part (a) 
may be applied to $\Ga^*$ yielding that 
$G(B_i, B_j)$ is transitive on $(X_{ji}\times X_{ij})\setminus A(\Ga)$
whenever $X_{ij}\ne\emptyset$.
Conversely supose that  $G(B_i, B_j)$ is transitive on 
$(X_{ji}\times X_{ij})\setminus A(\Ga)$ for some $i, j$ such 
that $X_{ij}\ne\emptyset$. 
Since $\Ga$ is $G$-symmetric, and since $\Ga^*$ is not null, 
it follows that $A(\Ga^*)\cap (B_i\times B_j) = (X_{ji}\times X_{ij})
\setminus A(\Ga)$ is non-empty. Thus there exists an arc 
$(x,y)$ of $\Ga^*$ with $x\in B_i$ and 
$y\in B_j$.  Let $(x',y')$ be arbitrary arc of $\Ga^*$, say
with $x'\in B_{i'}$ and $y'\in B_{j'}$, so also $X_{i'j'}\ne\emptyset$. 
By the definition of $\Ga^*$ there is an 
edge of $\Ga$ between $B_i$ and $B_j$, and similarly 
there is an edge of $\Ga$ between the parts $B_{i'}$ and $B_{j'}$.
Since $G$ is transitive on $A(\Ga)$, some element $g\in G$ maps 
$(B_{i'}, B_{j'})$ to $(B_i, B_j)$, and  
replacing $x',y'$ by their images under $g$ 
we may assume that $x'\in B_i$ and $y'\in B_j$.
Thus both $(x,y)$ and $(x',y')$ lie in $(X_{ji}\times X_{ij})\setminus A(\Ga)$, 
and transitivity of $G(B_i, B_j)$ on this set shows that  
$(x,y)$ and $(x',y')$ lie in the same $G$-orbit. 
Thus $\Ga^*$ is $G$-symmetric, proving (b).

Now $t\leq k$ (see Table~\ref{tab-param}), and if $t=k$ then
by its definition, $\Ga^*$ is a null graph. Suppose that $t<k$
and that $\Ga^*$ is $G$-symmetric. 
To verify the assertions about the parameters (see Table~\ref{tab-param}) 
we observe that $\Ga_\B=\Ga^*_\B$, 
and hence $b^*=b$. Also $v^*=|B_i|=v$ and $k^*=|X_{ij}|=k$ 
(for $X_{ij}\ne\emptyset$). From Definition~\ref{def-2.6},
we see that $t^*=k-t$.
Since $t<k$, for $x\in B_i$, say,
the set of parts of $\B$ containg vertices of $\Ga(x)$ is 
the same as that containing vertices of $\Ga^*(x)$, so $r^*=r$,
and hence, by \eqref{eq-1}, $s^*=r^*t^* = r(k-t)=rk-s$, completing the proof.
\qed

The only instance of the above construction that we shall need
arises when the pair $(\G, \B)$ has $t=1$ and $G(B_i,B_j)$ is
2-transitive on $X_{ij}$. In such a setting we obtain a 
natural correspondence between certain triples $(\G, G, \B)$ 
with $t=1, k\ge 3$ and their $*$-transforms $(\G^*, G, \B)$
with $t^*=k-t=k-1\geq 2$. We note that when $t=1$ and $k=|X_{ij}|=2$,
then $G(B_i,B_j)$ is automatically 2-transitive on $X_{ij}$,
so $\G^*$ is always $G$-symmetric, and is precisely the graph
called $\G^{opp}$ in~\cite[Section~6]{GP1}.

\begin{example}\label{ex-2.7}%
{\rm
The $*$-transform of the triple $(\G, G, \B)$ in Example~$\ref{ex-2.4}$ 
of valency $6$ (with $k=3, t=2$) is a triple $(\G^*, G, \B)$ of
valency $3$ with $t^*=1$ and $\G_\B^*=\G_\B=K_5$ (in fact 
$\G^* = 5\cdot K_4$). 
The $*$-transform of the triple $(\G, G, \B)$ of valency $6$ in 
Example~$\ref{ex-2.5}$ (with $k=3, t=2$) is a triple $(\G^*, 
G, \B)$ of valency $3$ with $k^*=3, t^*=1$, and $\G_\B^* =\G_\B = K_8$; 
vertices $(i,\b)$ and $(i',\b')$ are adjacent in $\G^*$ if and only if 
$\b=\b'$, so $\G^*=14\cdot K_4$.
}
\end{example}

\begin{example}\label{ex-2.8}%
{\rm
Example~$\ref{ex-2.5}$ generalises naturally to flags in
$\AG(d,q)$, or in $\PG(d,q)$, giving $G$-symmetric triples 
$(\G, G, \B)$ with $\G_\B$ a complete graph and with $G(B)^B$
$2$-transitive, for $B\in\B$. The resulting triples are 
$*$-transforms of simpler triples with $t^*=1$ and with  
the same quotient. (These examples show that the parameter $t$
can be arbitrarily large, as does the natural partition $\B$ 
in the complete multipartite graphs $K_{(b+1)[v]}$ with $b+1$ parts 
each of size $v$.)
}
\end{example}

\section{Imprimitive symmetric graphs with $t\ge 1$}\label{sec-3}%

Let $\G$ be an arbitrary $G$-symmetric, imprimitive graph with vertex set $V$.
If $x\in V$, then $\G_i (x)$ denotes the set of vertices at distance $i$ 
from $x$, for $i\ge 1$, and we usually write $\G (x)=\Ga_1(x)$.
Let $\B =\{ B=B_0, B_1, B_2, \ldots \}$ be a $G$-invariant partition of $V$.  If $B_i,
B_j$ are two parts joined by an edge of $\G$, then (since $G$ acts
transitively on directed edges of $\G$) each vertex $x\in B_i$ is joined to
either 0 or $t$ vertices of $B_j$, where $t\ge 1$ is independent of the choice
of the parts $B_i, B_j$.

In~\cite{GP1} we considered only the case $t=1$. In fact we assumed that
$G(x)$ acts primitively on $\G (x)$, whence either the valency $b$ of
the quotient graph $\G_\B$ is 1 and
$\G$ is bipartite, or $t=1$. However the approach developed in~\cite{GP1}
can often be used in situations where $t\ge 2$. In this section we run
through the basic ideas following~\cite[Sections~3 and~4]{GP1}. The proofs
are straightforward and are omitted.

By our notational convention, if $B\in\B$, then $\G_\B (B)$ denotes
the set of ``vertices'' in the graph $\G_\B$ which are adjacent to
the ``vertex'' $B$ in the graph $\G_\B$. 
That is to say, $\G_\B (B)$ denotes the set of
parts in $\B$ which are joined to $B$ by an edge of $\G$. 
Recall the definition of $\D(B)$ at the end of Subsection~\ref{sub:basic}.
The next result justifies the descriptions of the parameters in Table~\ref{tab-param}.

\begin{proposition}\label{prop-3.1}%
Let $\G$ be a $G$-symmetric graph whose vertex set $V$ admits a 
$G$-invariant partition $\B =\{B_0,
B_1, B_2,\ldots \}$. 
\begin{description}
\item[(a)] There exists a constant $t\ge 1$ such that, for each 
$x\in B$ and each $B_i\in\G_\B (B)$, the cardinality $|\G (x)\cap B_i|=0$ or $t$.
\item[(b)] $\G$ is regular of valency $s$ (say), and $t$ divides $s$.
Set $r:=s/t$; then each point of $\D (B)$ lies in $r$ blocks of $\D (B)$.
\item[(c)] Each pair of adjacent parts is joined by a constant number
$m$ (say) of edges of $\G$, and $t$ divides $m$. Set $k:=m/t$; then if
$B_i\in\G_\B (B)$, $k=|\{ x\in B\ :\ \G (x)\cap
B_i\ne\emptyset\} |$ is the number of points incident with each block of 
$\D (B)$.
\item[(d)] Each part $B\in\B$ has constant size $v$ (say), and the
quotient graph $\G_\B$ is regular of valency $b$ (say), where $vs=bm$
(so $vr=bk$).
\item[(e)] $\D (B)$ is a $1$-design with parameters $(v,b,r,k)$.
\item[(f)] $G(B)$ induces a group of automorphisms of $\D (B)$, which
is transitive on ``points'' (that is, on $B$), on ``blocks'' (that
is, on $\G_\B (B)$), and on ``flags'' (that is, on $\{(x,B_i)\ :\ x\in B,
B_i\in\G_\B(B)\ \mbox{and}\ \G(x)\cap B_i\ne\emptyset\}$).
\end{description}
\end{proposition}

\medskip\noindent
{\bf Corollary 3.1.1}\quad {\it Let $\G$ be a $G$-symmetric graph
and let  $\B =\{ B=B_0, B_1, B_2,\ldots \}$ be a $G$-invariant partition of the
vertex set. Suppose in addition that $G(B)^B$ is $2$-transitive.
\begin{description}
\item[(a)] Then each pair of ``points'' of $B$ lies in a constant
number $\l$ (say) of ``blocks'' of $\D (B)$. Hence either $\l =0$ 
 or $\D (B)$ is a $2$-design -- possibly with repeated
``blocks''.
\item[(b)] $\l =0$ if and only if $k=1$ (whence $t=k=m=1$, $b=vr$).
\item[(c)] Suppose that $\l\ge 1$ (whence $k\geq 2$) and that 
$\D (B)$ has repeated
``blocks''. Then each ``block'' is repeated the same number 
of times, say $\rho$, so $\rho\ge 2$ and $\rho$ divides $\l , r$, and $b$;
if we ignore repetitions then $\D (B)$ becomes a $2-(v,k,\l /\rho
)$ design with $b/\rho$ ``blocks'', so either
\begin{description}
	\item[(i)] $v=k,\ \rho =r=b=\l$; or
	\item[(ii)] $v>k$, whence $v\le b/\rho <b$ (by Fisher's inequality).
\end{description}
\item[(d)] Suppose that $\l\ge 1$ and that $\D (B)$ has no repeated
``blocks'' (so $\rho =1$ and $k\ge 2$). Then either $b=1$ or 
we must have $v>k$, so $v\le b$ (by Fisher's inequality). Hence 
\begin{description}
	\item[(i)] $b=1, v=k$, (so, if $\G_\B$ is connected, then  
	$\G =v\cdot K_2$ or $\G=K_{v,v}$); or
	\item[(ii)] $k<v<b$; or
	\item[(iii)] $k<v=b$, $\D (B)$ is a non-degenerate symmetric
$2$-design, and $G(B)$ acts $2$-transitively on both ``points'' and
``blocks'' of $\D (B)$. 
\end{description}
\end{description}
}
\medskip\noindent 
{\bf Corollary 3.1.2}\quad {\it Let $\G$ be a
$G$-symmetric graph and let  $\B =\{ B=B_0, B_1, B_2,\ldots \}$ be a
$G$-invariant partition of the vertex set. Suppose that $G(B)^B$ is
$2$-transitive, and that, in addition, the quotient graph 
$\G_\B=K_{b+1}$ is a complete graph. Then one of the following holds. 
\begin{description}
\item[(a)] $v<b$; 
\item[(b)] $v=b$, $\l =0$, $t=k=m=1=r$, so $\G = \binom{v+1}{2}
\cdot K_2$, and (as in {\rm \cite[Theorem $4.2$]{GP1}}) $V=\{ ij :\ 0\le i, j\le
v, i\ne j\}$, $B_i=\{ ij:\ 0\le j\le v, j\ne i\}$ ($0\le i\le v$), vertex
`$ij$' is joined only to `$ji$', 
and $G$ may be any $3$-transitive group on $\{ 0, 1,\ldots ,v\}$; 
\item[(c)] $v=k$ and $\D (B)$ has a single ``block'' of size $v$
repeated $b=r$ times;
\item[(d)] $v=b,\ 2\le k<v$, and $\D (B)$ is a symmetric $2$-design
with no repeated blocks, whence $G(B)$ is $2$-transitive on
``blocks'' of $\D (B)$ (that is, on $\G_\B (B)=\B \setminus\{ B\}$) as
well as on ``points'', so $G$ is $3$-transitive on $\B$.
\end{description}
}

\medskip If we wish to classify all such triples $(\G, G, \B)$ with $v\ge b$,
it remains to analyse the triples occurring in parts (c) and (d) of
Corollary~3.1.2.  In~\cite{GP2} we proved that when $t=1$ in case (c), then
either $v\le b+1$ or $(\G , G, \B)$ is uniquely determined; we also made
inroads into classifying the triples with $v=b+1$ or $v=b$.  Thus the following 
problem remains open. Its solution would complete the classification of all 
exceptional triples $(\Ga, G, \B)$ with $v\ge b$.

\medskip\noindent 
{\bf Problem 3.1.3}\quad
\emph{Classify all triples $(\Ga, G,\B)$ satisfying the conditions of 
Corollary $3.1.2$~(c) with $v\ge b$ and $t\ge 2$.}

\medskip

The rest of this paper is devoted to classifying the triples $(\G , G, \B )$ in part (d) of
Corollary~3.1.2: Section~\ref{sec-4} classifies the triples with $t=1$, while
Section~\ref{sec-5} classifies those with $t\ge 2$.

Note that when $k<v$ and $G(B)^B$ is 2-transitive, the 
group $G$ is faithful on $\B$ (see the first paragraph of 
the proof of Lemma~\ref{lem-3.2} below). Let $G_B$ denote the pointwise 
stabiliser of the part $B$. If
$\D (B)$ has no repeated blocks, then $G_B$ fixes each part in
$\G_\B(B)$, so if $\G_\B$ is a complete graph then $G_B=1$; hence 
$G(B)$ is faithful on $B$. Thus when classifying triples which arise
in part (d) of Corollary~3.1.2 we do not need to worry about
unfaithful  actions. (This is very different from case (c) of
Corollary~3.1.2; see~\cite{GP2}.) 

In part (d) of Corollary~3.1.2, $G^\B$ is 3-transitive and the subgroup $G(B)$
is 2-transitive both on $\B\setminus\{ B\}$ and on $B$ of the same degree
$v=b$.  The case $k=v-1$ is dealt with in Proposition~\ref{prop-4.1} (for $t=1$) and
Proposition~\ref{prop-5.1} (for $t\geq 2$).  If $2\le k\le v-2$, or equivalently
$2\le r\le b-2$, then a vertex $x\in B$ will be adjacent to vertices from $r$
parts of $\B\setminus\{ B\}$, and we note that $|\B\setminus\{ B\}| = b \geq r+2$. 
Consequently the stabiliser $G(x)$ of a point
$x\in B$ is not transitive on $\B\setminus\{ B\}$ (since $r<b$) and, although $G(x)$ 
has the same order as the stabiliser in $G_B$ of a part of $\B\setminus\{ B\}$,
$G(x)$ is not equal to the stabiliser of such a part $B'\in \B\setminus\{ B\}$ 
(for if it were then $G(x)$ would be transitive on $\B\setminus\{ B, B'\}$
and this is not the case since $r\le b-2$).  In particular the
actions of $G(B)$ on $\B\setminus\{ B\}$ and on $B$ are not equivalent.  The
next lemma identifies the possibilities for $G$ explicitly.

\begin{lemma}\label{lem-3.2}%
Let $\G$ be a
$G$-symmetric graph and let  $\B =\{ B=B_0, B_1, B_2,\ldots \}$ be a
$G$-invariant partition of the vertex set. Suppose that,  $G(B)^B$ is
$2$-transitive, that the quotient graph $\G_\B=K_{b+1}$ is a
complete graph, and that part (d) of Corollary~$3.1.2$ holds
with $k\le v-2$. Then $G$ is a
group of automorphisms of a $3$-design $\D$ with point set $\B$ such
that, for $\b$ a block of $\D$ containing $B$, $G (B,\b )$ is
transitive on $\b\setminus\{ B\}$ and on $\B\setminus \b$.  Moreover,
$(\D, G)$ are as in one of the lines $2$--$4$ of Table~$\ref{tab-main2}$. 
\end{lemma}

\proof 
First we show that $G$ acts faithfully on $\B$. Let
$K$ be the kernel of this action. If $K\ne 1$, then $K^B$ is a
nontrivial normal subgroup of the 2-transitive group $G(B)^B$,
so $K^B$ is transitive. Thus, for $x\in B$, $G(B)=KG(x)$. Hence
$G(x)^\B = G(B)^\B$ is transitive on $\B\setminus\{ B\}$,
so $r=b$ whence $k=v$ (since $vr=bk$), contrary to Corollary~3.1.2(d). 
Therefore $K=1$. 

It follows from the remarks immediately before Lemma~\ref{lem-3.2}
that $(G,\B )$ is a 3-transitive group  of degree
$b+1$ such that, for $B\in\B$, $G(B)$ has two 2-transitive
permutation representations of degree $b$ which are not
equivalent and which are such that a point stabiliser in one
of the representations is intransitive in the other representation.
Hence from the classification of finite 2-transitive groups
(see~\cite{C} and~\cite[Appendix 1]{L}) the socle of $G(B)$ 
is one of the following: 
$\PSL (d,q)$ with $d\ge 3$ and $b=(q^d-1)/(q-1)$; or $A_7$ with
$b=15$; or $HS$ with $b=176$; or $\PSL (2,11)$ with $b=11$.
A 3-transitive extension of a group $G(B)$ in the first family exists
if and only if either $q=2$ and $G=\AGL (d,2)$ as in 
Table~\ref{tab-main2}~line 4, or
$(d,q)=(3,4)$ and $G=M_{22}$ or $\Aut (M_{22})$ as in Table~\ref{tab-main2}~line 2. 
Similarly $A_7$ has a
unique transitive extension $G=2^4\cdot A_7$ as in Table~\ref{tab-main2}~line 4, 
and $\PSL (2,11)$ has a unique transitive extension $G=M_{11}$ as in
Table~\ref{tab-main2}~line 3. The group $HS$ has no transitive extension. 
In each case $G$ acts on a 3-design $\D$ with point set $\B$ as claimed, and
$G(B,\b )$ is transitive on both $\b\setminus\{ B\}$ and
$\B\setminus\b$, where $\b$ is a block of $\D$ containing $B$.\qed
\medskip

In each of the cases of Lemma~\ref{lem-3.2} we may therefore
identify $\B$ with the
point set of the 3-design $\D$. We finish this section with a lemma
which shows that the vertices of $\G$ may be identified with the
flags of $\D$ and that the symmetric 2-design $\D (B)$ is related to
$\D$ in one of two ways. 

\begin{definition}\label{def-designs}
{\rm
Let $\D$ be an $s$-design, where $s\geq 1$, and 
let $\calP$ denote the point set of $\D$. 
\begin{enumerate}
 \item[(a)] For a point $P$ of $\D$  the {\it derived design} $\D_P$ is the
2-design with point set $\calP\setminus\{ P\}$ and with blocks the sets
$\b\setminus\{ P\}$, where $\b$ is a block of $\D$ containing $P$.

\item[(b)] The {\it dual design} of $\D$ is the design whose points are the
blocks of $\D$, and whose blocks are the points of $\D$, with
incidence unchanged. The dual design of $\D$ is, in general, 
a 1-design, and is a 2-design if and only if $\D$ is a symmetric 2-design.

\item[(c)] The {\it complementary design} $\D^c$ of $\D$ has the same 
point set $\calP$ as $\D$, with blocks of $\D^c$ being the complements
(in $\calP$) of those of $\D$.
\end{enumerate}
 }
\end{definition}

\begin{lemma}\label{lem-3.3}%
Let $\G$ be a $G$-symmetric graph and let  
$\B =\{ B_0, B_1, B_2,\ldots \}$ be a
$G$-invariant partition of the vertex set. Suppose that, for $B\in\B$, $G(B)^B$ is
$2$-transitive, that the quotient graph $\G_\B=K_{b+1}$ is a complete
graph, and that
case (d) of Corollary~$3.1.2$ holds with $k\le v-2$. Then the
following also hold.
\begin{description}
\item[(a)] $(G,\B)$ is a $3$-transitive permutation group which is equivalent 
to the action of $G$ on the point set $\calP$ of a $3$-design
$\D$ as in Table~$\ref{tab-main2}$, lines $2$--$4$.

\item[(b)] The vertex set $V$ may be identified with the set of flags
of $\D$ in such a way that $G$ acts coordinate-wise; each part
$B_P$ of $\B$ is the set of flags with a fixed first coordinate $P$.

\item[(c)] For $B=B_P\in\B$ corresponding to a point $P\in\calP$,
there are just two possibilities for the non-degenerate $2$-design 
$\D(B_P)$: either
	\begin{description}
	\item[(i)] $\D (B_P)$ is the dual design of the derived design
$\D_P$ (that is, the point set $B_P$ of $\D (B_P)$ corresponds to the set 
of blocks of $\D_P$, and the block set $\G_\B(B_P)=\B\setminus\{ B_P\}$ of
$\D (B_P)$ corresponds to the set of points of $\D_P$; moreover if $\b$ 
is a block of $\D$ incident with $P$, 
then the point $(P,\b )$ in $B_P$ is incident with a block
$B_{P'}$ of $\D (B_P)$ if and only if $P'\in\b$); or
	\item[(ii)] $\D (B_P)$ is the complement of the 
dual design of the derived design $\D_P$: that is, a point $(P,\b )$ of $\D
(B_P)$ (which corresponds to a block $\b$ of $\D$ incident with $P$) 
is adjacent to a
block $B_{P'}$ of $\D (B_P)$ (which is a point $P'$ of $\D_P$) if and
only if $P'\not\in\b$.
	\end{description}
\item[(d)] Let $x=(P,\b )\in B_P\in\B$, and let
$B_{P'}\in\B\setminus\{ B_P\}$ be such that $\G (x)\cap
B_{P'}\ne\emptyset$.
	\begin{description}
	\item[(i)] In case (c)(i) (that is, if $P'\in\b$), the
$G(x,B_{P'})$-orbits in $B_{P'}$ are $\{ (P',\b )\}$, $\{ (P',\b')
:\ P,P'\in\b'\ \mathrm{and}\ \b'\ne\b\}$, and $\{ (P',\b') :\ P'\in\b',
P\not\in\b'\}$. The orbit lengths are as follows: $1, 4, 16$  in line $2$ of Table~$\ref{tab-main2}$; 
and $1, 4, 6$ in line $3$ of Table~$\ref{tab-main2}$; and $1, 2^{d-1}-2, 2^{d-1}$  
in line $4$ of Table~$\ref{tab-main2}$. 
	\item[(ii)] In case (c)(ii) (that is, if $P'\not\in\b$), the
$G(x,B_{P'})$-orbits in $B_{P'}$, in  line $2$ of Table~$\ref{tab-main2}$, 
are $\{ (P',\b') :\ P,P'\in\b'\}$ of
length $5$, $\{ (P',\b') :\ P'\in\b', \b'\cap\b =\emptyset\}$ of
length $6$, and $\{ (P',\b') :\ P\not\in\b', P'\in\b',
|\b'\cap \b |=2\}$ of length $10$.
In  lines $3$ and $4$ of Table~$\ref{tab-main2}$, the
$G(x,B_{P'})$-orbits in $B_{P'}$ are as follows:
$\{ (P',\bar{\b} )\}$ (where $\bar{\b}=\calP\setminus\b$),
$\{ (P',\b') :\ P\not\in\b', P'\in\b', \b'\ne\bar{\b}\}$, and $\{
(P',\b') :\ P,P'\in\b'\}$, of lengths $1, 5, 5$ in line $3$, and of lengths $1,
2^{d-1}-1, 2^{d-1}-1$ in line $4$.
	\end{description}
\end{description}
\end{lemma}

\proof
By our comments before Lemma~\ref{lem-3.2}, and by 
Lemma~\ref{lem-3.2}
itself, we may identify $\B$ with the point set $\calP$ of one of the
3-designs $\D$ of Table~$\ref{tab-main2}$, lines $2, 3, 4$. 
Let $B=B_P\in\B$ correspond
to the point $P\in\calP$. As we noted before the statement of 
Lemma~\ref{lem-3.2},
the 2-transitive actions of $G(B)$ on $\B\setminus\{ B\}$ and on $B$
are not equivalent and it follows that the action of
$G(B)$ on $B$ is isomorphic to its action on the blocks of $\D_P$. It
follows that we may label each point $x\in B$ as $(P,\b )$, where
$\b\setminus\{ P\}$ is the block of $\D_P$ to which it corresponds,
in such a way that $G(B)$ acts naturally on the second coordinates.
Thus we may label elements of the set $V$ by flags $(P,\b )$ of $\D$
in such a way that $G$ acts coordinate-wise and each part of $\B$
consists of the set of flags with a fixed first coordinate. Thus parts (a) and (b) are proved.

Let $x=(P,\b )\in B$. By Lemma~\ref{lem-3.2}, $G(x)=G(P,\b )$ is transitive on
both $\b\setminus\{ P\}$ and $\calP\setminus \b$, and since the
$G$-actions on $\B$ and $\calP$ are equivalent, 
$G(x)$ has two orbits in $\B\setminus\{ B\}$, namely $\B_1:=\{ B_{P'}
:\ P'\in\b , P'\ne P\}$ and $\B_2:=\{ B_{P'} ;\
P'\in\calP\setminus\b\}$. It follows, since $G(x)$ is transitive on
$\G (x)$, that $\G (x)$ meets either each block in $\B_1$ or each
block in $\B_2$ (but not both of these since $r\le b-2$). Now points of
$B_P$ correspond to blocks of $\D_P$; and the block set of $\D
(B_P)$, namely $\B\setminus\{ B_P\}$, is the point set of $\D_P$. If
$\G (x)$ meets the parts of $\B_1$, then a point $(P,\b )$ and a
block $B_{P'}$ of $\D (B_P)$ are incident in $\D (B_P)$ if and only
if $P'\in\b\setminus\{ P\}$, and therefore $\D (B_P)$ is the
dual design of $\D_P$. On the other hand if $\G (x)$ meets the parts
of $\B_2$ then $\D (B_P)$ is the dual design of $\D_P$ with incidence
reversed - that is, $\D(B_P)$ is the complement of the dual of $\D_P$. 

To prove part (d) we make a careful examination of the cases.
Suppose first that $\D (B_P)$ is the dual design of $\D_P$ as in case
(c)~(i). Let $x=(P,\b )\in B_P$ and $P'\in\b\setminus\{ P\}$. Then
in all the cases of Lemma~\ref{lem-3.2}, $G(x,B_{P'})=G(P,P',\b )$ has three
orbits in $B_{P'}$, namely the fixed point $(P',\b )$, the set of all
$(P',\b')$ where $\b'$ contains both $P$ and $P'$, and the set of all
$(P',\b')$ where $\b'$ contains $P'$ but not $P$. The orbit lengths
are as stated in (d)~(i). Now suppose that $\D (B_P)$ is the dual
design of $\D_P$ with incidence reversed, as in case (c)~(ii).
This time consider $x=(P,\b )\in B_P$ and 
$P'\not\in\b$. In  lines $3$ and $4$ of Table~$\ref{tab-main2}$,
$\bar{\b}=\calP\setminus\b$ is also a block of $\D$ and $P'\in\bar{\b
}$. Moreover $G(x,B_{P'})=G(P,P',\b )$ also fixes $(P',\bar{\b})\in
B_{P'}$ and so $G(x,B_{P'})=G(B_P, (P',\bar{\b}))$; it can be
checked that the $G(P,P',\b )$-orbits in $B_{P'}$ are $\{
(P',\bar{\b})\}$, $\{ (P',{\b'}) ;\ P,P'\in\b'\}$, and $\{ (P',{\b'})
;\ P\not\in\b', P'\in\b', \b'\ne\bar{\b}\}$. The orbit lengths are as
stated in (d)~(ii).

In  line $2$ of Table~$\ref{tab-main2}$, $G(B_P,B_{P'})=G(P,P')$ has orbits of
lengths 5, 16 in $B_P$ and in $B_{P'}$, and $x=(P,\b )$ lies in a
$G(P,P')$-orbit of length 16 (since $P'\not\in\b$). Thus (since 5 and
16 are coprime) $G(x,B_{P'})=G(P,P',\b )=A_5$ is still transitive on
the $G(P,P')$-orbit of length 5 in $B_{P'}$. We claim that this is
the set of vertices $(P',\b')$ such that $\b'$ contains $P, P'$ and one
point of $\b\setminus\{ P\}$. By~\cite[p.39]{Atlas} there are 60 blocks of
$\D$ which meet
$\b$ in two points; each such block contains 4 points of
$\calP\setminus\b$ and consequently (since $G(\b )$ is transitive on
$\calP\setminus\b$) the point $P'$ of $\calP\setminus\b$ lies
in 15 of these blocks. Since each triple of points of $\D$ lies in a
unique block of $\D$, there are exactly 5 blocks of $\D$ containing
$P ,P'$ and one point of $\b\setminus\{ P\}$. Thus $P'$ lies in
exactly 10 blocks of $\D$ which do not contain $P$ and which meet $\b$
in 2 points. Moreover, from~\cite[p.39]{Atlas}, we see that there are 16
blocks of $\D$ which are disjoint from $\b$ and, as $G(\b )$ is
transitive on
$\calP\setminus\b$, it follows that $P'$ lies in exactly 6 blocks
disjoint from $\b$. It follows that the three subsets listed in
(d)~(ii) in this case are all invariant under $G(P,P',\b )
\cong A_5$. From our discussion above, we can now conclude 
that the set of size 5 is an orbit, and, since
$G(P,P',\b )$ is transitive on the unordered pairs from
$\b\setminus\{ P\}$, the set of size 10 is also an orbit. From the
character table for $M_{22}$ in~\cite[p. 40]{Atlas}, we see that elements of
$G(P,P',\b )$ of order 3 fix exactly 5 blocks of $\D$. Since they
must fix two blocks of the 5 containing $P$ and $P'$, and they must
fix one block of the 10 containing $P'$ but not $P$ and meeting $\b$
in two points, they can fix at most two blocks of the 6 containing
$P'$ and disjoint from $\b$. It follows that they fix none of the
latter blocks and that the set of size 6 is also a $G(P,P',\b
)$-orbit. \qed
 
\section{$\D (B)$ a symmetric 2-design with $t=1$}\label{sec-4}%

We assume througthout this section that $(\G ,G,\B )$ is a triple
arising in case (d) of Corollary~3.1.2 with $t=1$. We begin by
looking briefly at the case $k=2<v$. This is not strictly necessary
but gives some insight into the general case where $k=v-1$. We then classify triples $(\G ,G,\B )$ with
$k=v-1$, before analysing the non-degenerate case $3\le k\le v-2$. 

\medskip\noindent\emph{Case $k=2<v$:}\quad
Here each of the $\binom{v}{2}$ pairs of ``points''
in $B$ lies in $\l$ ``blocks'' of $\D (B)$, so $\binom{v}{2}\cdot \l =
b=v$. Hence $v=3, \l =1, r=k=2$, so $\B =\{ B_0, B_1, B_2, B_3\}$ and
each vertex $x\in B_i$ is joined to one vertex in each of two of the
other three parts, but to no vertex of the third part, say $B_{i'}$.
Thus each vertex $x\in V$ receives a natural label $ii'$ ($0\le i\ne
i'\le 3$), where $x\in B_i$ and $\G (x)\cap B_{i'}=\emptyset$. In
particular no vertex whose label has second coordinate $i'$ is joined
to any vertex whose label has $i'$ as its first
coordinate. The
$G$-partition
$\B$ is determined by this labelling; $G(B)=S_3$, so $G=S_4$ acting
coordinate-wise on labels, as in Theorem~\ref{main}~(b). There are just two possibilities for
$\G$: either 
\begin{description}
\item[(i)] any two adjacent vertices have the same second coordinate
(and different first coordinates), so $\G = 4\cdot K_3$; or
\item[(ii)] any two adjacent vertices have labels involving four
different coordinates, so $\G = 3\cdot C_4$.
\end{description}

The case $k=2$ is instructive in that $k=2$ implies $k=v-1$ (since $v=3$ in (i) and (ii)),
and both the approach used when $k=2$, and the outcomes (i) and (ii) 
constitute a useful `first approximation' to the general case $k=v-1$.
Possibility (i) ($\G = 4\cdot K_3$) occurs in line 2 of Table~\ref{tab-main1}.  
We interpret the graph $\G = 3\cdot C_4$ in (ii) above in two different 
ways.  The first interpretation is to
observe that, since $b+1=4$, the parts $B_i$ could have been indexed
with the points of the projective line $\PG (1,3)$; the vertices of
$\G$ would then have been labelled by ordered pairs of points from
$\PG (1,3)$, and $\infty 0$ would then be joined to both $12$ and
$21$, whence $\G$ is visibly the cross ratio graph $\CR (3;1)$ (as
defined in Subsection~\ref{sub-cr}); this exhibits $3\cdot C_4$ as the
first member of the family of graphs in line 4 of Table~\ref{tab-main1}. 
The second interpretation is to view the coordinates $\{ 0,1,2,3\}$ in
the vertex labels as the four points of an affine geometry $\AG
(2,2)$ of dimension $d=2$, with $ij$ joined to $i'j'$ in $\G$ if and
only if $i,j,i',j'$ are pairwise distinct and lie in a single
2-dimensional subspace! This exhibits $3\cdot C_4 = (2^2-1)
\cdot K_{2^{2-1}[2]}$ as the first member of the family of graphs in 
Table~\ref{tab-main1} line $5$ with $\D'$ as 
in Table~\ref{tab-main2} line $1$.

\begin{proposition}\label{prop-4.1}%
Let $\G$ be a $G$-symmetric graph
with vertex set $V$. Let $\B =\{B_0, B_1, B_2, \ldots 
\}$ be a $G$-partition of $V$ with $|\B|=b+1$ and each $|B_i|=v$,
and with quotient 
$\G_\B =K_{b+1}$. Suppose that $v=b$ and that, for $B\in\B$, $G(B)^B$ is
$2$-transitive. Suppose
further that $t=1$ and that either $k=2<v$ or $k=v-1\ge 2$. Then the
vertex set $V$ may be taken to be $\{ ij :\ 0\le i, j\le v, i\ne
j\}$, where $x=ij$ if and only if $x\in B_i$ and $\G (x)\cap
B_j=\emptyset$. The $G$-invariant partition $\B$ is determined by this
labelling and $\G , G, \B$ are as in Theorem~$\ref{main}$~(b) satisfying one of 
Table~$\ref{tab-main1}$ line $2$, $4$, or $5$ (with $\D'$ as 
in Table~\ref{tab-main2} line $1$). 
%
\end{proposition}

\proof
We have already seen that this result is true when
$k=2$, so we may assume that $k=v-1\ge 3$. Since $v=b$, we have $r=k=v-1$ by
Proposition~\ref{prop-3.1}~(d), so each vertex $x\in B_i\in\B$ is joined 
to one vertex in each of $r=v-1$ parts $B_{i'}\ne B_i$. Moreover, if $B_{i'}\ne B_i$
and
$\G (x)\cap B_{i'}=\emptyset$, then $x$ is the only vertex of $B_i$
not joined to $B_{i'}$ (since $k=v-1$). Thus each vertex $x\in V$
receives a natural label $ii'$ as claimed with $G$ acting
coordinate-wise. In particular, the $G$-partition is determined by
the labelling, and no vertex whose label has second coordinate $j$
can be joined to a vertex whose label has first coordinate $j$. This
leaves just two possible cases; we consider each in turn.

Suppose two adjacent vertices have labels with the same second
coordinate; that is, suppose that some vertex $il\in B_i$ is joined
to $jl\in B_j$, where $i\ne j$. Let $B:=B_l=\{ ll' : \ 0\le l'\le v,
l'\ne l\}$. Then 
$G(B)^B$ is 2-transitive, and the actions of $G(B)$ on $X=\{ il :\
0\le i\le v, i\ne l\}$ and on $B$ are equivalent. Hence the
induced subgraph $\la X\ra = K_v$, so
$\G$ is the disjoint union of $v+1$ copies of $K_v$, with $G^\B$
3-transitive, as in Table~\ref{tab-main1} line $2$. 
Conversely, once the vertices have been labelled, the
graph $\G$ in which vertices are joined if they have the same second
coordinate has
$\Aut (\G ,\B )=S_{v+1}$, so $G$ may be any 3-transitive subgroup of
$S_{v+1}$. 

There remains the possibility that each pair of adjacent vertices $x,
y$ have labels involving four different symbols, say  
$x=01,\ y=cd$. Let $B=B_0 = \{ 0i\ :\ 1\le i\le v\}$. Each
coordinate symbol ``$i$''  corresponds to a part $B_i\in\B$.
As in the previous paragraph, $G(B)$ is
2-transitive on the set $X=\{ i0 : \ 1\le i\le v\}$ and hence $G$ is
3-transitive on $\B$, that is, on the coordinate symbols $\{ 0,1,\ldots
,v\}$. Also $G(x)=G(01)$ is transitive on the remaining $v-1$
vertices in $B$, namely $\{ 0i :\ 2\le i\le v\}$. Let $w=0c$,
for some $c>1$. Then $G(x,w)=G(01c)$ fixes $c0$, so leaves $B_c$ 
invariant; hence $G(x,w)$ fixes $\G (x)\cap B_c=\{ y =cd\}$,
say, so fixes $d$. Thus $G$ must be a
3-transitive group on $\{ 0,1,\ldots ,v\}$ in which the stabiliser
$G(01c)$ of three points fixes a fourth point $d$. Checking the list
of 3-transitive groups (for example in~\cite{C, L}) we see that there are
just two possible cases: either 
\begin{description}
  \item[(i)] $\PSL (2,q)\le G\le\PGaL (2,q)$ for some prime power $q\ge
3$, or
  \item[(ii)] $G=\AGL (d,2)$ for some $d\ge 2$, or $G=Z_2^4\cdot A_7<\AGL
(4,2)$ with $d=4$.
\end{description}
The graphs arising in case (i) are classified in \cite[Theorem 4.1]{cr}.
Since $k\geq 3$ and since adjacent vertices have labels involving four 
different symbols, the only examples are cross-ratio graphs $\CR
(q;d,s)$ or $\TCR(q;d,s)$ as defined in Definition~\ref{def-2.1}.
Since $t=1$ it follows (see Remark~\ref{rem-cr}~(d)) that $s=s(d)$.
Thus in this case $\G, G$ satisfy  line $4$ of Table~$\ref{tab-main1}$.

Finally suppose that $G=\AGL (d,2)$ or $Z_2^4\cdot A_7<\AGL (4,2)$.
The vertices of $\G$ are labelled by ordered pairs of distinct points
of the affine geometry $\AG (d,2)$. Moreover we showed above that in
this case, if $01$ is joined to $cd$, then $G(01c)$ fixes $d$; thus
$01$ is joined to $cd$ if and only if $0, 1, c, d$ are the four
points of a 2-dimensional subspace. Hence  line $5$ of 
Table~\ref{tab-main1} holds with $\D'$ as in  line $1$ of Table~\ref{tab-main2}.
Note that, even when $G=Z_2^4\cdot A_7$, $G(01c)$ has only four fixed points in $\AG
(d,2)$, so the possibilities for $G$ are also as listed. 
(Here $\G = \ColPairs(\AG_2(d,2))\cong (2^d-1)\cdot K_{2^{d-1}[2]}$,
and each 2-dimensional subspace $\{ 0, 1, c, d\}$ 
in $\AG(d,2)$ gives rise to three cycles of length 4 in $\G$, 
for example, $\la 01, cd, 10, dc\ra = C_4$.) 
\qed

\medskip
We end this section by considering the remaining case $3\le k\le v-2$, $t=1$.

\begin{proposition}\label{prop-4.2}%
Let $\G$ be a $G$-symmetric graph with
vertex set $V$. Let $\B =\{B_0, B_1, B_2, \ldots \}$
be a $G$-invariant partition of $V$ with $|\B|=b+1$ and each $|B_i|=v$,
 and with quotient $\G_\B 
=K_{b+1}$. Suppose that $v=b$ and that, for $B\in\B$, 
$G(B)^B$ is $2$-transitive. Suppose further that $3\le k\le v-2$, and $t=1$.
Then we may take $V$ to be the set of flags $P\b$ in a
$3-(v+1,k+1,\l )$ design $\D$, where
$\G , G, \B, \D$ are as in Theorem~$\ref{main}$~(c), and moreover satisfy one of 
Table~$\ref{tab-main3}$ line $1$ (with $\D$ as 
in Table~$\ref{tab-main2}$ line $2, 3$ or $4$), or line $2$ (with $\D$ as 
in Table~$\ref{tab-main2}$ line $3$ or $4$).
\end{proposition}

\proof
Now $\D (B)$ is a symmetric design with $3\le k\le
v-2$, and part (d) of Corollary~3.1.2 holds. It follows from
Lemma~\ref{lem-3.3} that $V$ may be identified with the flags of a 3-design
$\D$ such that each part $B\in\B$ is the set of flags on a certain point of $\D$; 
moreover $\D$  is as in one of lines 2--4 of Table~\ref{tab-main2}.
Thus each vertex $x$ is identified with a flag $P\b$, where 
$P$ is a point of $\D$ and $\b$ is a block of $\D$ incident with $P$,
and hence the first assertions of Theorem~\ref{main}~(c) hold.
Since $t=1$, it follows that, if $x= P\b$ is adjacent to
$P'\b'$, then $P\ne P'$ and $G(x, B_{P'})=G(P,P',\b )$ must fix $P'\b'$ and
hence must fix $\b'$. By Lemma~\ref{lem-3.3}~(d), it follows that either 
\begin{description}
  \item[(i)] $P\b$ is adjacent to $P'\b'$ if and only if $\b
=\b'$ and $\b$ contains both $P$ and $P'$; or
  \item[(ii)] $P\b$ is adjacent to $P'\b'$ if and only if $\b$
and $\b'$ are disjoint.
\end{description}
In case (i) the graph $\G$ is the disjoint union of $b'$ complete
graphs $K_{k'}$, where $b'$ is the number of blocks of $\D$ and
$k'$ is the number of points of $\D$ incident with a given
block. It follows from Lemmas~\ref{lem-3.2} and~\ref{lem-3.3} 
that $\G, G, \B, \D$ satisfy line 1 of Table~\ref{tab-main1}, 
with $\D$ as in  one of lines 2--4 of Table~\ref{tab-main2}. 
On the other hand in case~(ii) above, the graph $\G$
is the disjoint union of $b'/2$ complete bipartite graphs
$K_{k',k'}$, and by Lemmas~\ref{lem-3.2} and~\ref{lem-3.3}, 
line 2 of Table~\ref{tab-main1} holds, 
with $\D$ as in  one of lines 3 or 4 of Table~\ref{tab-main2}.\qed

\section{$\D (B)$ a symmetric 2-design with $t\ge 2$}\label{sec-5}%

We assume throughout this section that $(\G , G,\B )$ is a triple
arising in case (d) of Corollary~3.1.2 with $t\ge 2$. We begin by
looking briefly at the cases $k=2 (=r)$ and $k=3 (=r)$. Both the
analysis and the examples that arise emphasise the importance of the
vertex labelling; we therefore make this labelling explicit at the
outset.  When $t\ge 2$, the `degenerate case' $k=v-1$
behaves differently for different groups $G$; hence our classification 
of the triples $(\G ,G,\B )$ with $k=v-1$ (Proposition~\ref{prop-5.1})
proceeds by considering each family of groups separately.
Finally we deal with the non-degenerate case $3\le k\le v-2$
(Proposition~\ref{prop-5.2}) where we have available to us the results of
Lemmas~\ref{lem-3.2} and~\ref{lem-3.3}.

\medskip\noindent
\emph{Case $k=2$:}\quad Here each
pair of points of $B$ lies in a unique block of $\D (B)$ (and,
since Corollary~3.1.2 (d) holds, 
there are no repeated blocks), so $\binom{v}{2}=b=v$. Hence
$v=3$, $k=v-1$, $\B =\{ B_0, B_1, B_2, B_3\}$, $G(B)=S_3$, and
$G=S_4$. As at the beginning of Section~\ref{sec-4}, each vertex $x$ of $\G$
receives a unique label $x=ii'$, where $x\in B_i$, $i\ne i'$, and $\G
(x)\cap B_{i'}=\emptyset$. Let $x=01$. Then $G(x)=G(01)=Z_2$ must
act transitively on $\G (x)$, so $|\G (x)\cap B_2|=t\le 1$. Hence
$k=2$ does not occur when $t\ge 2$. 

\medskip\noindent
\emph{Case $k=3$:}\quad 
Here each pair of points in $B$ lies in $\l\ge 1$
blocks of $\D (B)$, so $\binom{v}{2}\l = b\binom{k}{2}=v\cdot
3$. Hence (since $k=3<v$) either $v=4, \l =2$, or $v=7, \l =1$.
Suppose first that  $v=4, \l =2$; then $k=v-1$, $\D (B)$ is the
complete $2-(4,3,2)$ design, $G(B)=A_4$ or $S_4$, $G=A_5$ or $S_5$.
Since $k=v-1$, each vertex $x$ of $\G$ receives a unique label
$x=ii'$ ($0\le i, i'\le 4$, $i\ne i'$), where $x\in B_i$, 
and $\G (x)\cap B_{i'}=\emptyset$. Let $x=01$. Then $G(x,
B_2)$ fixes $21$ and $20$, so $t\ge 2$ implies that $\G (x)\cap B_2=\{
23, 24\}$ and that $G=S_5$. So $t=2$, and $\G =\Del^*$ is obtained as
the ``$*$-transform'' of a pair  $(\Del ,\B )$ with $t=1$ from
Proposition~\ref{prop-4.1}. By uniqueness of $\Ga$, the graph 
$\Del = (4+1)\cdot K_4$ (recall the definition of $\Del^*$ from 
Definition~\ref{def-2.6}), so $\Ga, G$ satisfy Theorem~\ref{main}~(b),
namely line 3 of Table~\ref{tab-main1}. Note that $\G$ also occurs as
the cross-ratio graph $\CR(4;d,1)$, where $d\in\GF(4)\setminus\{0,1\}$
 in line 4 of Table~\ref{tab-main1}.

Suppose next that $v=7, \l =1$. Since
$k=3$, $\D (B)$ is $\PG (2,2)$ and $G(B)=\PSL (3,2)$. Since $G$ is a
transitive extension of the 2-transitive group ``$G(B)$ acting on the
lines of $\PG (2,2)$'', we have $G=\AGL (3,2)$. The eight parts of
$\B$ are the points of $\D =\AG (3,2)$. Each line of $\D (B)$ corresponds to a ``point''
$B_i\in\B\setminus\{ B\}$; and each point $x$ of $\D (B)$ is
determined by the three lines of $\D (B)$ which contain it, that is, the
three points $B_i\in\B\setminus\{ B\}$ such that $\G (x)\cap
B_i\ne\emptyset$. Since $k=3=r$, these three points $B_i$ in
$\B\setminus\{ B\}$, together with $B$, form a hyperplane of $\D
=\AG_2 (3,2)$. Thus each vertex $x$ is naturally labelled by a flag
$B\b$ of $\D =\AG_2(3,2)$, where $x\in B$, and $\b$ is the 
hyperplane of $\D$ which contains $B$ and which gives rise to the line
$\b\setminus\{ B\}$ which corresponds to $x$ in the derived design 
$\D_B \cong \PG(2,2)$. Thus $x=B\b$ is joined to $x'=B'\b'$ only
if $B\ne B'$, $B\in\b'$ and $B'\in\b$. If the labels of adjacent
vertices $x, x'$ could have the same second coordinate $\b =\b'$, then $G(x,
B')$ would fix $x'$; but $G(x,B')$ must act transitively on $\G
(x)\cap B'$, so $t=|\G (x)\cap B'|=1$. Hence when $t\ge 2$, 
$x=(B,\b )$ is joined to
$x'=(B',\b')$ only if $\b\ne\b'$; so since $G(x,B')$ is transitive on
the two hyperplanes which contain $B$ and $B'$ and which are not
equal to $\b$, we have $t=2$ and adjacency in $\G$ is determined.
This time $\G =\Del^*$ is the ``$*$-transform'' of a pair $(\Del ,\B )$
with $t=1$ from Proposition~\ref{prop-4.2}.  By uniqueness of $\Ga$, the graph  
$\Del = (2^4-2)\cdot K_{2^2}$, so $\Ga, G$ satisfy Theorem~\ref{main}~(c),
namely line 1 of Table~\ref{tab-main3} with $\D$ in line 4 of Table~\ref{tab-main2}. 

\medskip
These two examples illustrate what we expect to find in general:
$k=3, v=4, \l =2$ is an instance of the `degenerate'  case $k=v-1$
(Proposition~\ref{prop-5.1}); $k=3, v=7, \l =1$ is an instance of the
non-degenerate  case $3\le k\le v-2$ (Proposition~\ref{prop-5.2}).
Many of the triples $(\G ,G,\B )$ arising in Propositions~\ref{prop-5.1} 
and~\ref{prop-5.2} are the $*$-transforms of triples with $t=1$ 
from Section~\ref{sec-4}. In addition we shall obtain  
cross-ratio graphs and a few sporadic examples. 

We begin by considering the `degenerate'  case $k=v-1$. Here the
design $\D (B)$ is the complete $2-(v,v-1,v-2)$ design.

\begin{proposition}\label{prop-5.1}%
Let $\G$ be a $G$-symmetric graph
with vertex set $V$. Let $\B =\{B_0, B_1, B_2, \ldots 
\}$ be a $G$-partition of $V$ with $|\B|=b+1$ and each $|B_i|=v$, 
and with quotient 
$\G_\B =K_{b+1}$. Suppose that $v=b$ and that, for $B\in\B$, $G(B)^B$ is
$2$-transitive. Suppose
further that $t\ge 2$ and that either $k=2<v$, or $k=v-1\ge 2$. Then
$k>2$ and the vertex set $V$ may be taken to be $\{ ij :\ 0\le i, j\le
v, i\ne j\}$, where $x=ij$ if and only if $x\in B_i$ and $\G (x)\cap
B_j=\emptyset$. 
Further $\G , G, \B$ are as in Theorem~$\ref{main}$~(b), and satisfy one of 
Table~$\ref{tab-main1}$ line $3$ or $4$, or line $5$ or $6$ (with $\D'$ as 
in Table~$\ref{tab-main2}$ line $1, 2$ or $3$).

\end{proposition}

\proof
We have already seen that $k\ne 2$ and that, when $k=3$, the 
conclusion holds (with $v=4, G=S_5$). 
Thus we may assume that $v-1=k\ge 4$. Since
$k=v-1$, $\D (B)$ is the complete $2-(v,v-1,v-2)$ design.
In particular, for any $B_j\in\B\setminus\{ B_i\}$, there is a unique
$x\in B_i$ such that $\G (x)\cap B_j=\emptyset$. Conversely, for each
$x\in B_i$, there is a unique $B_j\in\B\setminus\{ B_i\}$ such that $\G
(x)\cap B_j=\emptyset$. Thus each $x$ receives a natural label $ij$
as claimed with $G$ acting coordinate-wise. In particular the
$G$-invariant partition $\B$ is determined by the labelling, and no vertex
whose label has second coordinate $j$ can be joined to a vertex whose
label has first coordinate $j$.

For $x=ij\in B_i$, the stabiliser $G(x)$ fixes $B_i$ and $B_j$
setwise: since $G$ acts coordinate-wise on labels, $G(B_i,B_j)=
G(ij) = G(x)$. Thus the actions of $G(B_i)$ on $\B\setminus\{ B_i\}$
and on $B_i$ are equivalent. In particular (since
$G(B_i)^{B_i}$ is 2-transitive) $G$ is 3-transitive on $\B$, and 
(since $G$ acts coordinate-wise on labels) $G$ acts faithfully on
$\B$. Now $G(x)$ is transitive on $\B\setminus\{ B_i, B_j\}$ and on $\G
(x)$, so $\G (x)$ consists of $t$ vertices from each part $C\in
\B\setminus\{ B_i, B_j\}$. Moreover, $\G (x)\cap C$ is an orbit of
$G(x,C)$ of length $t\ge 2$. Now $G(x,C)$ is the stabiliser in $G$ of
the three parts $B_i, B_j$ and $C$. Consequently $G(x,C)^C$ is the
stabiliser in the 2-transitive group $G(C)^C$ of two vertices, say
$u,u'$, and $\G (x)\cap C$ is a $G(x,C)$-orbit in $C\setminus\{
u,u'\}$ of length $t\geq 2$. In particular, $G(B_i, B_j, C)^C=G(x,C)^C\ne 1$.
 If $C=B_l$ then we may assume that $u, u'$ are the vertices $li, lj$, respectively. 

Suppose first that $G(x,C)$ is transitive on $C\setminus\{ u,u'\}$,
that is, that $G$ is 4-transitive on $\{ 0,1,\ldots v\}$ (so $G$ is
one of the groups listed in line 3 of Table~\ref{tab-main1}). 
Then $\G (x)\cap C =
C\setminus\{ u,u'\}$, that is, $ij$ is adjacent to $ll'$ if and only
if $l'$ is distinct from $i, j, l$. So $\G$ is as in line 3 of Table~\ref{tab-main1}, that is,
$\G$ is the $*$-transform of the graph in in line 2 of Table~\ref{tab-main1}.

Thus we may suppose that $G$ is 3-transitive, but not 4-transitive, on
$\B$, and also that the stabiliser of 3 points in this action is
nontrivial (since $G(B_i, B_j, C)^\B\cong G(B_i, B_j, C)\ne 1$).  
Then $G$ is one of the following groups.
	\begin{description}
	\item[(i)] $\PSL (2,q)\le G\le\PGaL (2,q)$ on the $v+1=q+1$ points of
the projective line $\PG (1,q)$, such that $G$ is
3-transitive, and $G(\infty 01)\ne 1$ (so $q\ge 8$ since we are assuming that
$v=q>4$); or
	\item[(ii)] $G=\AGL (d,2)$ with $v+1=2^d$ for some $d \ge 3$, 
or $v+1=2^4$ and $G=Z_2^4\cdot A_7<\AGL (4,2)$; or
	\item[(iii)] $G=M_{22}$ or $\Aut (M_{22})$ with $v+1=22$; or
	\item[(iv)] $G=M_{11}$ with $v+1=12$.
	\end{description}
 Suppose first that $\PSL (2,q)\le G\le\PGaL (2,q)$, where $q=p^n$
for some prime $p$ and $n\ge 2$. Then the vertices of
$\G$ are labelled by ordered pairs of distinct points from the
projective line $\PG (1,q)$. If the vertex $x$ is labelled $\infty 0$,
then $\G(x)\cap B_1$ must be a single $G(\infty 01)$-orbit of
length $t\ge 2$. Thus if $y\in \G(x)\cap B_1$, then $y=1d$, where 
$d\ne \infty, 0, 1$. It follows from \cite[Theorem 4.1]{cr} that
$\G$ is a cross-ratio graph $\CR(q;d,s)$ or $\TCR(q;d,s)$ for some divisor
$s$ of $s(d)$, where the subfield of $\GF(q)$ generated by $d$ has order
$p^{s(d)}$. Moreover, by \cite[Theorems 3.4 and 3.7]{cr}, $t=s(d)/s\geq 2$,
and hence $\G, G$ satisfy line 4 of Table~\ref{tab-main1}.

Next suppose that $G=\AGL (d,2)$ ($d\ge 3$) or $Z_2^4\cdot A_7<\AGL
(4,2)$. The vertices of $\G$ are labelled by ordered pairs of distinct
points of the affine geometry $\AG (d,2)$. If $01$ is joined to $cd$,
then $G(01c)$ does not fix $d$ (since $t\ge 2$); hence $d\not\in\{ 0,1\}$
and $d$ is not equal to the
fourth point $c'$ of the affine plane spanned by $0, 1, c$. Now 
$G(01c)$ is transitive on $\AG (d,2)\setminus\{ 0,1,c,c'\}$. This is
easily checked if $G=\AGL (d,2)$, while if $G=Z_2^4\cdot A_7$, then
$G(01c)=A_4$ and each nonidentity element of this subgroup has
exactly three fixed points in $\AG (4,2)\setminus\{ 0\}$ 
(see~\cite[p. 10]{Atlas}) so it must act regularly on the 12 points of $\AG
(4,2)\setminus\{ 0,1,c,c'\}$. Hence $01$ is joined to $cd$ if and
only if $0, 1, c, d$ are distinct points and span an affine space of
dimension 3.  Thus $\G, G$ are as in line 6 of 
Table~\ref{tab-main1} (with $\D'$ as in line 1 of 
Table~\ref{tab-main2}). 

Now we treat the cases where $G=M_{22}$ or $\Aut (M_{22})$ ($v+1=22$),
or $M_{11}$ ($v+1=12$). Here the vertices of $\G$ are labelled
by ordered pairs of distinct points of the design $\D$, where $\D$ 
is  the 3-(22,6,1) Steiner system or the unique 3-(12,6,2) design
respectively. Suppose that $x=ij\in B_i$ is joined to $y=i'j'\in
B_{i'}$. Then  $G(iji')$ has two fixed points in $B_{i'}$, namely $i'i$
and $i'j$, and two nontrivial orbits, namely the set $C_1$ of
vertices $i'l$ where $i,j,i',l$ lie in some block of $\D$, and
the set $C_2$ of vertices $i'l$ where $i,j,i',l$ are contained
in no block of $\D$. Thus in each case we
have exactly two graphs, corresponding to $\G (x)\cap B_{i'}$ equal
to $C_1$ or $C_2$ respectively, and  line 5 or 6 of 
Table~\ref{tab-main1} holds (with $\D'$ as in line 2 or 3 of 
Table~\ref{tab-main2}). 
\qed

\medskip
Finally we treat the ``non-degenerate'' case where $3\le k\le v-2$.

\begin{proposition}\label{prop-5.2}%
Let $\G$ be a $G$-symmetric graph
with vertex set $V$. Let $\B =\{B_0, B_1, B_2, \ldots 
\}$ be a $G$-partition of $V$ with $|\B|=b+1$ and each $|B_i|=v$, 
and with quotient 
$\G_\B =K_{b+1}$. Suppose that $v=b$ and that, for $B\in\B$, $G(B)^B$ is
$2$-transitive. Suppose
further that $3\le k\le v-2$, and $t\ge 2$. Then
we may take $V$ to be the set of flags in a
$3-(v+1,k+1,\l )$ design $\D$, where
$\G , G, \B, \D$ are as in Theorem~$\ref{main}$~(c), and moreover satisfy one of 
Table~$\ref{tab-main3}$ line $3$ (with $\D$ as 
in Table~$\ref{tab-main2}$ line $2, 3$ or $4$), or line $4$ (with $\D$ as 
in Table~$\ref{tab-main2}$ line $3$ or $4$), or  line $5$ or $6$ (with $\D$ as 
in Table~$\ref{tab-main2}$ line $2$).
\end{proposition}

\proof
The 2-design $\D (B)$ is symmetric 
with $3\le k\le v-2$, and part (d) of Corollary~3.1.2 holds. It
follows from Lemma~\ref{lem-3.3} that $V$ may be identified with the flags of a
3-design $\D$ as claimed, where $\D$ is as in one of lines 2-4 of Table~\ref{tab-main2}.
Moreover, since $t\ge 2$, each vertex $x=P\b$ is adjacent to
all the points of some nontrivial orbit (of length $t$) of $G(x,
B_{P'})=G(P,P',\b )$ in $B_{P'}$. By Lemma~\ref{lem-3.3}~(d), it follows
that either 
	\begin{description}
	\item[(i)] $\D (B_P)$ is the dual design of the derived design 
$\D_P$, so that, for $P'\ne P$, we have  $\G (B_{P'})\cap B_P =\{ P\b' :\
P'\in\b'\}$; or
	\item[(ii)] $\D (B_P)$ is the dual design of the derived design 
$\D_P$ with incidence reversed so that, for $P'\ne P$,  
$\G (B_{P'})\cap B_P =\{ P\b' :\ P'\not\in\b'\}$.
	\end{description}
In case (i), $x=P\b$ is adjacent to a vertex of $B_{P'}$ if and
only if $P'\in\b$. Thus if $x=P\b$ is adjacent to
$P'\b'$, then we must have $P\in\b'$. Moreover,by
Lemma~\ref{lem-3.3}~(d)~(i), $G(x, B_{P'})$ is transitive on the set of all
points $P'\b'$ such that $\b'$ contains $P$, and hence in this case,
$P\b$ is adjacent to $P'\b'$ if and only if $P\ne P'$,
$\b\ne\b'$, and $\b , \b'$ contain both $P$ and $P'$. Thus
$\G$ is the $*$-transform of the
corresponding graph in Proposition~\ref{prop-4.2}, namely of a graph from line 1 of Table~\ref{tab-main3},
and hence line 3 of Table~\ref{tab-main3} holds for $\Ga$, and 
the values of $t$ can be read off from the orbit lengths in 
Lemma~\ref{lem-3.3}~(d)~(i).

In case (ii), $x=P\b$ is adjacent to a vertex of $B_{P'}$ if and
only if $P'\not\in\b$. Thus if $x=P\b$ is adjacent to
$P'\b'$, then we must have $P\not\in\b'$. Suppose first that $\D$ is
in line 3 or 4 of Table~\ref{tab-main2}, that is $\D$ is either  the unique
3-(12,6,2) design or the affine
geometry of points and hyperplanes of $\AG (d,2)$, respectively. 
Then by Lemma~\ref{lem-3.3}~(d)~(ii), $G(x,
B_{P'})$ is transitive on the set of all vertices $P'\b'$ such that
$P\not\in\b'$ and $\b'$ is not the complement $\bar{\b}$ of $\b$, and
$G(x, B_{P'})$ fixes $P'\bar{\b}$. Since $t>1$ it follows that 
$P\b$ is adjacent to $P'\b'$ if and only if $P\ne P'$,
$P\not\in\b'$, $P'\not\in\b$, and $\b'$ is not the complement
$\bar{\b}$ of $\b$. This graph $\G$ is the $*$-transform of the
corresponding graph in Proposition~\ref{prop-4.2}, namely of a graph from line 2 of Table~\ref{tab-main3},
and hence line 4 of Table~\ref{tab-main3} holds for $\Ga$, and 
the values of $t$ can be read off from the orbit lengths in 
Lemma~\ref{lem-3.3}~(d)~(ii). 

Finally suppose that $\D$ is the
3-(22,6,1) Steiner system as in line 2 of Table~\ref{tab-main2}. Then by
Lemma~\ref{lem-3.3}~(d)~(ii), $G(x, B_{P'})$ has two orbits on the set of all
vertices $P'\b'$ such that $P\not\in\b'$, namely the set of those
for which $\b\cap\b'=\emptyset$ (of length $t=6$) and the set of those
for which $|\b\cap\b'|=2$ (of length $t=10$). Thus there are two
possible rules for incidence, giving two graphs $\G$ in this case. 
For the first graph, $P\b$ is adjacent to $P'\b'$ if and only if
$P\ne P'$, $P\not\in\b'$, $P'\not\in\b$, and $\b\cap\b'=\emptyset$,
while for the second graph, $P\b$ is adjacent to $P'\b'$ if and
only if $P\ne P'$, $P\not\in\b'$, $P'\not\in\b$, and
$|\b\cap\b'|=2$. Thus we have line 5 or 6 of Table~\ref{tab-main3} respectively. \qed

\medskip
The results of this section complete the proofs of Theorem~\ref{main},
noting that the examples in Table~\ref{tab-main1} line 1 arise from 
Corollary 3.1.2 (b), and that the assertions of Theorem~\ref{main}~(b)
hold for these examples.

\end{document}